\documentclass[a4paper,10pt]{article}
\usepackage{amsmath,amsthm,amssymb,amsfonts,epsfig,epstopdf,titling,url,array,enumerate,cite,mathrsfs,upgreek,authblk,scrextend, graphicx, float}
\usepackage[bindingoffset=0.2in,left=1in,right=1in,top=1in,bottom=1in,footskip=.25in]{geometry}
\theoremstyle{plain}
\newtheorem{thm}{Theorem}[section]
\providecommand{\keywords}[1]{\begin{addmargin}[28pt]{28pt}\noindent\textbf{Keywords:} #1 \end{addmargin}}
\newtheorem{lma}[thm]{Lemma}
\newtheorem{cor}[thm]{Corollary}
\newtheorem{ppn}[thm]{Proposition}
\theoremstyle{definition}
\newtheorem{dfn}[thm]{Definition}
\newtheorem{eg}[thm]{Example}
\newtheorem{rem}[thm]{\textit{Remark}}
\providecommand{\ams}[1]{\begin{addmargin}[28pt]{28pt}\noindent\textbf{Mathematics Subject Classification:} #1\end{addmargin}}

\title{\textbf{A new approach to generalize  metric functions}}
\author[1]{Abhishikta Das}
\author[2]{Anirban Kundu}
\author[3,*]{T. Bag}
\affil[1,2,3]{Department of Mathematics, Siksha-Bhavana,	\authorcr	Visva-Bharati, Santiniketan-731235, West-Bengal, India  
	\authorcr 	E-mail\textsuperscript{1}: abhishikta.math@gmail.com
	\authorcr  E-mail\textsuperscript{2}: anirbankundu92@gmail.com 
	\authorcr E-mail\textsuperscript{3,*}: tarapadavb@gmail.com}

\date{}

\begin{document}
	\maketitle
	\begin{abstract}
		\noindent
		S-metric and b-metric spaces  are metrizable, but it is still quite impossible to get an explicit form of the concerned metric function.
		To overcome this, the notion of $\phi$-metric is developed by making a suitable modification in triangle inequality and its
		properties  are pretty similar to metric function. It is shown that one can easily  construct a $\phi$-metric from existing generalized 
		distance functions  like S-metric, b-metric, etc. and those  are $\phi$-metrizable. The convergence of sequence on those metric spaces
		is identical to the respective induced $\phi$-metric spaces. So, unlike metrics, concerned $\phi$-metric can be easily constructed 
		and  $\phi$-metric functions may play the role of metric functions substantially. Also, the structure  of $\phi$-metric spaces is studied 
		and some fixed point theorems are established. 
	\end{abstract}
	%%%%%%%%%%%%%%%
	%%
	%%
	%%
	%%%%%%%%%%%%
	%%%%%%%%%%%%%%%%%%%%%%%%%%
	\ams{47H10, 54H25.}
	%%%%%%%%%%%%%
	%%%
	%%
	%%%
	%%%
	%%%
	%%%%%%%%%%%%%%
	%%%%%%%%%%%%%%%%%%%%%%%%%%
	\section{Introduction } 
	In modern mathematics metric spaces and topological spaces are two widely used concepts. Metric spaces are considered 
	as a particular case of topological spaces. The notion of  metric was developed by Frechet\cite{[13]} and later Hausdorff \cite{[14]}
	presented axiomatically. The proposed three axioms   are very fundamental and geometrically understandable. The properties 
	of metric spaces are easier to check than  general topological spaces. Because of this reason, metrizability is an interesting topic 
	for topological spaces. Unfortunately, there are spaces which are not metrizable.   So researchers keep on developing functions which are
	more general than metric spaces which are not metrizable. In 1963, Gahler\cite{[7]} introduced 2-metric and later n-metric in general. 
	These spaces are not metrizable in general.  In another process of generalization, Dhage introduced D-metric\cite{[17]}. But it was defective  
	structure. Later  Mustafa and Sims introduced G-metric spaces \cite{[18]},  and Sedghi introduced S-metric spaces\cite{[1]}. On the other hand,
	another kind of space, called b-metric space (\cite{[3]}, \cite{[2]}) was introduced. 2-metric, G-metric, and S-metric were defined to generalize the concept 
	of distance between two points to the distance between three points. All of these three concepts generalize the three axioms of
	metric functions in their own way. But b-metric is a generalization of metric with straight modification in the axiom, called
	‘triangle inequality’. \\ 
	In 2013, Chaipuniya and Kumam introduced the notion of g-3ps \cite{[70]} and claimed that it is the general structure of the distance between
	three points. They proved the    G-metric and S-metric spaces are also g-3ps and a g-3ps is not metrizable in general. But it
	has been proved that   b-metric, G-metric and S-metric spaces are metrizable. So we think that instead of g-3ps, which is 
	much more general in nature, a general structure closer to those distance functions may be a better alternative.  \\
	With this idea in mind, we introduce the notion of $\phi$-metric, a generalized distance function, and show the b-metric and S-metric 
	spaces
	are particular examples only. As we have mentioned earlier, the above-mentioned spaces are metrizable. In fact, we show that 
	$\phi$-metric spaces are metrizable as well. The line of proof for metrizability shows that one can not construct
	a metric for b-metric and S-metric spaces easily.
	But  $\phi$-metric can be constructed easily and serves well in the absence
	of an easily constructible metric.  Thus, one does not need to find the metric for an S-metric or b-metric. Instead, $\phi$-metric 
	will be enough to study the topological structure. \\
	%%%%%%%%%%%%%%%%%%%%%%
	The organization of this article is as follows. Section 2 provides some preliminary results. In Section 3,  $\phi$-metric spaces are defined, illustrated  by examples and some elementary topological properties are studied. In Section 4, results on compactness, completeness, and totally boundedness are discussed.  Section 5 
	consists of some fixed point theorems in $\phi$-metric spaces. The straightforward proofs are omitted.
	%%
	%%
	%%%%%
	\section{ Preliminaries}
	In this section, we recollect some preliminary results which are used in this paper.
	%%%%%%
	\begin{dfn}\label{2metric} \cite{[7]}
		Let $X$ be a non-empty set. Then $(X, D)$ is called a 2-metric space if the function $ D : X \times X \times X \rightarrow \mathbb{R}$, named 
		2-metric satisfies the following conditions: 
		\begin{enumerate}[(i)]
			\item For every $a, b \in X$ with $a \neq b$ there exists $c \in X$ such that $D(a, b, c) \neq 0$;
			\item If at least two of $a, b, c \in X$ are the same, then $D(a, b, c) =0$;
			\item For all $a, b, c \in X$, 
			$D(a, b,c) = D(a, c,b) = D(b, c,a) = D(b, a,c) = D(c, a,b) = D(c, b,a)$;
			\item The rectangle inequality: for all $a, b, c, d \in X$, 
			$D(a, b,c) \leq D(a, b,d) + D(b, c,d) + D(c, a,d)$.
		\end{enumerate}
	\end{dfn}
	\begin{dfn}\label{bmetric}
		\cite{[3]} Let $X$ be a nonempty set and  $d : X \times X \rightarrow \mathbb{R}_{\geq 0}$ be a function,  
		for all $x, y, z \in X$ which satisfies the following conditions:  
		\begin{enumerate}[(b1)]
			\item  $d(x, y) = 0$ if and only if $x = y$;
			\item  $d(x, y) = d(y, x)$;
			\item $d(x, z) \leq 2[d(x, y) + d(y, z)]$.
		\end{enumerate}
		Then $ d $ is called a b-metric and the pair $( X, d )$ is called a b-metric space.
	\end{dfn}
	Later in 1998, Czerwik\cite{[2]} modified this notion of b-metric replacing $2$ by a constant $ K\geq 1 $ in the condition (b3) of Definition \ref{2metric}. Khamsi and Hussian \cite{[8]} took a step further and considered the constant $ K>0 $ and they named it as metric-type space.
	Another generalization of b-metric is Strong b-metric  space, which was introduced by Kirk and Shahzad \cite{[4]}.
	%
	%
	%
	%%%
	\begin{dfn}
		\cite{[4]} Let $X$ be a nonempty set and $K\geq 1$ be a given real number. A function $d : X \times X \rightarrow \mathbb{R}_{\geq 0}$ 
		is called a strong b-metric if for all $ x, y, z \in X $ it satisfies the followings:  
		\begin{enumerate}[(b1)]
			\item  $d(x, y) = 0$ if and only if $x = y$;
			\item  $d(x, y) = d(y, x)$;
			\item $d(x, z) \leq Kd(x, y) + d(y, z)$.
		\end{enumerate}
		Then $(X, d, K)$ is called a strong b-metric space.
	\end{dfn}
	%%%%
	\begin{dfn}
		\cite{[1]} Let $ X $ be a nonempty set. A function $S : X \times X\times X \rightarrow \mathbb{R} $ is called an S-metric 
		if it  satisfies the following properties:  
		\begin{enumerate}[(S1)]
			\item  $S(x, y,z) \geq 0,$ for all $x, y, z \in X $; 
			\item  $S(x, y,z)=0~$ if and only if $x = y=z$;
			\item $S(x, y,z) \leq S(x, x,w) + S(y,y, w)+S(z,z,w)$, for all $x, y, z, w \in X$.
		\end{enumerate}
		The pair $(X, S)$ is called an S-metric space.
	\end{dfn}
	%%%%
	%%%%
	\begin{ppn}\cite{[9]} Let $(X, d )$ be a b-metric space.
		\begin{enumerate}[(i)]
			\item A subset $A$ of  $X$  is called open if for any $a \in A,$  there exist $ t >0$ such that $B(a, t) \subset A$, where
			$$ B(a, t) = \{ y \in X : d( a, y) < t \}. $$
			\item If $\tau$ is the collection of all open balls of $(X, d)$, then $\tau $ defines a topology on $X$.
		\end{enumerate}
	\end{ppn}
	%%%%%%%%%%%%%%%%%%%%%%%%
	%%%%%%%%%%%%%%%%%%%%%%%%%%%%%%
	\begin{ppn}\cite{[4]}
		Every open ball $B(a, r) = \{ x \in X : d(a, x) < r \}$ in  a strong b-metric space $(X,d,K)$ is open.
	\end{ppn}
	%%%%%%%%%%%%%%%%%%%%%%%%%%%%%%%%%%%
	%%%%%%%%%%%%%%%%%%%%
	\begin{dfn} \cite{[1]}
		Let $(X, S)$ be an S-metric space. Then for any $x \in X $ and $t>0$, open ball and closed ball are defined by 
		$$ B_S (x, t)= \lbrace {y \in X : S( y,y,x) < t}\}~ \text{and}  ~~ B_S [x, t]= \lbrace {y \in X : S( y,y,x) \leq t}\} $$ respectively.
	\end{dfn}
	%%%%%%%%%%%%%%%%%%%%%%%%%%%%
	So far, we have discussed some generalized metric spaces.  In Section 3, we introduce  $\phi$-metric which is a generalization of many 
	established distance functions. Later, we involve ourselves to study the metrizability and topological properties of the discussed generalized
	spaces including $\phi$-metric spaces.  In this aspect, some topological definitions and results on topological spaces and  other 
	generalized metric spaces are below.    
	%%%%%%%%%%%%%%
	%%%%
	\begin{dfn} \cite{[10]}
		Let $X$ be a topological space and $ B =\{ B_s:s \in S \} $ be a family of subsets of $X$. Then 
		\begin{enumerate}[(i)]
			\item $B$ is called locally finite if, for each $x \in X$, there exists a neighborhood $U$ of $x$ such that the set
			$\{ s \in S: A_s\cap U\neq \phi \}$ is finite.
			\item $B$ is called discrete if, for each $x \in X$, there exists a neighborhood $U$ of $x$ such that the set
			$B_s\cap U\neq \phi$ for at most one $s\in S$.
			\item $B$ is called $\sigma$-locally finite if $B =\underset{i\in \mathbb{N}}{\cup}B_i$, where every $B_i$ is locally finite.
			\item $B$ is called $\sigma$-discrete if $B =\underset{i\in \mathbb{N}}{\cup}B_i$, where every $B_i$ is discrete.
			\item $B$ is called a cover of $X$ if $\underset{s\in S}{\cup}B_s=X$.
			\item A cover $A=\{A_i:i \in I\}$ of subsets of $X$ is called a refinement of the cover $B_i$ if for each $i \in I$, there 
			exists $s \in S$ such that $A_i\subset B_s$.
		\end{enumerate}	
	\end{dfn}
	%%%%%%%%%%
	%%%%%%%%%%%%%%%%%%%%%%%%%%%
	\begin{dfn}\cite{[10]}
		Let $(X,\tau)$ be a topological space. Then $X$ is said to be a 
		\begin{enumerate}[(i)]
			\item regular space if for any closed subsets $ A \subset X $ and for $ x \in X \setminus A $, there exist two disjoint open sets $U$ and $V$ 
			containing $A$ and $x$ respectively.
			\item normal space if for any two disjoint closed subsets $A$ and $B$ of $X$, there exist two disjoint open sets $U$ and $V$ 
			containing $A$ and $B$ respectively. 
		\end{enumerate}
	\end{dfn} 
	%%%%%%%%%%%%
	%%%%%%%%%%%%%%%%%%%%%%%%%
	\begin{dfn} (\cite{[11]},\cite{[12]}) Let $(X,\tau)$ be a topological space.
		\begin{enumerate}[(a)]
			\item A subset $U$ of $X$ is called sequentially open if each sequence $\{x_n\}\subset X$ converging to a point $x\in U$ then there
			exists $N\in \mathbb{N}$ such that $x_n\in U$, for all $n \geq N$.
			\item A subset $U$ of $X$ is called sequentially closed if no sequence in $U$ converges to a point, not in $U$.
			\item $X$ is called semi-metrizable if there exists a function $d :X\times X\rightarrow [0, \infty)$ such that for all $x, y\in X$,
			\begin{enumerate}[(i)]
				\item $d(x, y) =0$ if and only if $x =y$;
				\item $d(x, y) =d(y, x)$;
				\item $x \in \overline{A} $ if and only if $d(x, A) =inf\{d(x, y) :y\in A\} =0$ for any subset $A$  of $X $.
			\end{enumerate}
			\item $(X,\tau)$ is said to be metrizable if $ \exists$ a metric on $X$ whose topology is same as the topology $\tau$. 
		\end{enumerate}
	\end{dfn}
	%%%%%%%%%%%%%%%%
	%%%%%%%%%%%%%%%%%%%%%%%%%
	\begin{thm}\cite{[10]} \label{thm-metric-stone}
		(The Stone Theorem) Every open cover of a metrizable space has an open refinement which is both locally finite and $\sigma$-discrete. 	
	\end{thm}
	%%%%%%%%%%%%%%%
	%%%%%%%%%%%%%%%%%%%%%%%%%%
	\begin{thm}\cite{[15]} \label{thm-bing}
		(The Bing Metrization Theorem) A topological space is metrizable if and only if it is regular and has a $\sigma$-discrete base. 	
	\end{thm}
	%%%%%%%%%%%%
	%%%%%%%%%%%%%%%%%%%%%%%
	Recently, N. V. Dung et al.\cite{[9]} prove the following results about semi-metrizability, metrizability, and some other properties of b-metric spaces. 
	%%%%%%%%%%%
	%%%%%%%%%%%%%%%%%%%
	\begin{thm}\cite{[9]}
		Every b-metric space $(X, d)$ is a semi-metrizable space.	
	\end{thm}
	%%%%%%%%%%%%%%
	%%%%%%%%%%%%%%%%%%%%%%%%%
	\begin{thm}\cite{[9]} \label{thm1.2}	Let $ (X, d) $ be a b-metric space and $d$ is continuous in one variable. Then 
		\begin{enumerate}[(i)]
			\item $X$ is regular.
			\item  Every open cover of $X$ has an open refinement which is both locally finite and $\sigma$-discrete.
			\item $X$ has a $\sigma$-discrete base.
			\item $X$ is metrizable.
		\end{enumerate} 	
	\end{thm}
	%%%%%
	An important Corollary of Theorem \ref{thm1.2} is given below.
	%%%%%%%%%% 
	\begin{cor}\cite{[9]} \label{S is metrzble}
		Every S-metric space is metrizable. 	
	\end{cor} 
	%%%%%%%%%%%%%%%%%
	%%%%%%%%%%%%%%%%%%%%
	%%%%%%%%%%%
	%%%%%%%(Main Results)%%%%%%%%%%%%%%%
	%%%%%%%%%%%%%%%%%%%%%%%%%%%%%%%%%%%%%%%%%%%%%%
	\section{Introduction to $\phi$-metric}
	%%%%%%%%%%
	Recall that b-metric and S-metric spaces are metrizable (see Theorem \ref{thm1.2}, Corollary \ref{S is metrzble}). However, the role of induced 
	metrics in such spaces is implicit in nature.  Thus it is difficult to construct a metric for a given $ S $-metric or $ b $-metric so that 
	topology remains the same. \\
	In this connection, we  introduce $\phi$-metric in a new approach that helps to study the metrizability of $S$-metric spaces, b-metric spaces, etc. 
	%%%
	%%%
	%%%
	%%
	%%
	%%
	%%%%%%%%%%dfn of $\phi$- metric space%%%%%%%%%%%%%%%%%%%
	%%%%%%%%%%%%%%%%%%%%%%%%%%%%%%%%%%	
	\begin{dfn}
		Let $X$ be a nonempty set. A $\phi$-metric is a function $d_\phi:{{X \times X} \rightarrow}{\mathbb{R}_{\geq 0}}$  holding the following 
		conditions:
		\begin{enumerate}[($d_\phi1$)] \label{dfn1}
			\item $d_\phi(x,y)=0$ if and only if  $x=y$;
			\item $d_\phi(x,y)= d_\phi(y,x)$;
			\item $d_\phi(x,y) \leq d_\phi(x,z) + d_\phi(z,y)+\phi(x,y,z) $;
		\end{enumerate}
		for all $ x, y, z \in X$, where $ \phi : X \times X \times X \rightarrow \mathbb{R}_{\geq0} $ is a function satisfying 
		\begin{enumerate}[($\phi 1$)]
			\item  $\phi(x,y,z)=0 $ if $x= z$ or $y=z$;
			\item $\phi(x,y,z)=\phi(y,x,z)$;
			\item  For all  $ \epsilon>0 $, there exists  $  \delta>0$ such that $\phi(x,y,z)<\epsilon $ whenever $d_\phi(x,z)<\delta$ or $d_\phi(y,z)<\delta$;
		\end{enumerate}
		for all $ x, y, z \in X$. \\
		A set $ X $ together with the function $ d_\phi $,   $(X,d_\phi)$ is called a $\phi$-metric space. 
	\end{dfn}		
	%%%%%%%%%%%%%%%%%%%%%%%%%%%%%% EXAMPLE %%%%%%%%%%%%%%%%%%%%%%%%
	%%%%%%%%%%%%%%%%%%%%%%%%%%%%%%%%%%%%%%%%%%%%%%%%%%%%%
	%%%%%%%%%%%%%%%%%%%%%%%%%%%%%%%%%%%%
	\begin{eg}\label{ex1}
		Let $ (X,d) $ be a metric space and define a function on $X$ by $ d_\phi(x,y)=(d(x,y))^2 $,  for all  $ x,y\in X$. Then clearly
		$ d_\phi(x,y)=0 $ if and only if $ x=y$ and $ d_\phi(x,y)=  d_\phi(y,x),$ for all $ x,y\in X $. Now for any $ x,y, z\in X$, 
		\begin{align*}
			& d(x,y)  \leq d(x,z) + d(z,y)\\
			or ~~ & (d(x,y))^2 \leq (d(x,z))^2 + (d(z,y) )^2 + 2d(x,z)d(y,z)
		\end{align*}
		This implies $ d_\phi(x,y) \leq d_\phi(x,z) + d_\phi(z,y) + \phi(x,y,z) $
		where $\phi(x,y,z)=2\sqrt{d_\phi(x,z)d_\phi(z,y)},$ for all $ x,y,z\in X $.
		Hence $ (X,d_\phi) $ is a $\phi$-metric space.
	\end{eg}
	%%%
	\begin{rem}
		Every metric space is a $\phi$-metric space but not conversely. \\
		To justify, in Example \ref{ex1}, take $ X=l_p,~0<p<1 $ then the distance function $ d_\phi$ in $l_p,~0<p<1 $ is not a metric but it
		satisfies our new notion of distance. 	
	\end{rem}
	%%%%%%%%%%%%%%%%%
	\begin{eg}\label{ex2} 
		Let $(X,d_\phi)$ be a $\phi$-metric space and define a function on $X$ by 
		$D_{\phi_1}(x,y)=(d_\phi (x,y))^2, $ 
		for all $ x,y\in X$. 
		Then $ (X,D_{\phi_1}) $ is a  $\phi$-metric space.\\
		For clarification, proof of  inequality is enough. For any $x,y,z\in X$,
		\begin{align*}
			&  d_\phi(x,y)  \leq d_\phi(x,z) + d_\phi(z,y) + \phi(x,y,z) \\
			or ~~ & (d_\phi(x,y))^2 \leq (d_\phi(x,z))^2 + (d_\phi(z,y) )^2 + 2d_\phi(x,z)d_\phi(y,z) + (\phi(x,y,z))^2 + \\
			& ~~~~~~~~~~~~~~ ~~~~~~~~ 2 ( d_\phi(x,z) + d_\phi(z,y) ) \phi(x,y,z) 
		\end{align*}
		which implies $  D_{\phi_1}(x,y) \leq D_{\phi_1}(x,z) + D_{\phi_1}(z,y) + \phi_1(x,y,z) $	where $\phi_1(x,y,z)= 2 d_\phi(x,z )d_\phi(z,y)+ \phi^2(x,y,z)+ 2\phi(x,y,z)[d_\phi(x,z)+d_\phi(z,y)],$ 
		for all $ x,y,z\in X $.\\
		Moreover if $(X,d_{\phi_i}),~i=1,2\cdots n$ be $\phi$-metric spaces, then  $(X,D_\phi) $ is a $\phi$-metric space where 
		$D_\phi(x,y)=\prod_{i=1} ^n d_{\phi_i}(x,y),$ 
		for all $ x,y\in X$.\\
		But in case of infinite product of $\phi$-metrics that is  $D_\phi(x,y)= \underset{n\rightarrow\infty}{lim} (d_\phi(x,y))^n, $ 
		for all $ x,y\in X$ is a $\phi$-metric only if $d_\phi$  is the discrete metric.
	\end{eg}
	%%%%%%%%%%
	%%%%%%%%%%%%%%%%%%%%%%%%%%%
	%%%%%%%%%%%%%%%%%%%%%%%%
	\begin{eg}\label{ex3}
		Let $ (X,S)$ be a S-metric space.  Define 
		$ d_\phi(x,y)=S(x,x,y),$ for all $ x,y\in X$ and take  $\phi(x,y,z)=d_\phi(x,z)+d_\phi(y,z),$ 
		for all $ x,y,z\in X $. \\
		Then the function $\phi$ satisfies the first and second conditions. For the third condition,  take $\alpha >0$. 
		If $d_\phi(x,z)< \frac{\alpha}{2}$ and $d_\phi(y,z)< \frac{\alpha}{2}$, then $\phi(x,y,z)<\alpha$. Thus  
		for all $ \alpha >0, $ there exists $ \beta>0$ such that $ \phi(x,y,z)<\alpha$ whenever $d_\phi(x,z)<\beta $ 
		and $d_\phi(y,z)< \beta$ where $\beta= \frac{\alpha}{2}$. Hence $(X,d_\phi)$ is $\phi$-metric space.    		
	\end{eg}
	%%%%%%%%%%%%%%%%%%%%
	\begin{eg}\label{ex4}
		Let $ (X,B)$ be a b-metric space with constant coefficient $K(>1)$. For all $ x,y\in X $, define 
		$ d_\phi(x,y)=B(x,y),$ 
		and take  $\phi(x,y,z)=(K-1)[d_\phi(x,z)+ d_\phi(y,z)],$ 
		for all $x,y,z\in X $.\\
		To verify the third condition for $\phi$ function, take $\alpha >0$. Now if $d_\phi(x,z)< \frac{\alpha}{2K}$ and $d_\phi(y,z)< \frac{\alpha}{2K}$, then $\phi(x,y,z)< \dfrac{(K-1)\alpha}{K}< \alpha$. Thus 
		for all $\alpha >0, $ there exists  $ \beta>0$ such that $ \phi(x,y,z)<\alpha$ whenever $d_\phi(x,z)<\beta $ and $d_\phi(y,z)< \beta$ where $\beta=  \frac{\alpha}{2K}$. 
		Hence $(X,d_\phi)$ is $\phi$-metric space.	
	\end{eg}
	%%%%%%%%%%%%%%%%
	%%%%%%%%%%
	\begin{eg}\label{ex6}
		Let $ (X,B)$ be a strong b-metric space with constant coefficient $K(>1)$. Define $ d_\phi(x,y)=B(x,y), $ for all $ x,y\in X $,  and take $\phi(x,y,z)= (K-1)[ d_\phi(x,z)+ d_\phi(y,z)],$ for all $ x,y,z\in X $. 
		Then $(X,d_\phi)$ is $\phi$-metric space.	
	\end{eg}
	%%%%%%%
	\begin{rem}
		We call the $\phi$-metrics defined in the Example \ref{ex3}, Example \ref{ex4}, and Example \ref{ex6} as $\phi$-metric induced by S-metric, b-metric,  and strong b-metric respectively. So it is clear that one can easily construct a  $\phi$-metric from those generalized distance functions. 
	\end{rem}
	%%%%%%%%%%%%%%%%%%%%%%%%%%%%%%%%%%%
	%%%%%%%%%%%%%%%%%%%%%%%%%%%%%%%%%%%%%%%%%%
	%%%%%%%%%%%%%%%%%%%%%%%%% TOPOLOGY %%%%%%%%%%%%%%%%%%%%%%%%%%%%%
	To study the topological structure of  $\phi$-metric spaces, we define open and closed balls as given below. 
	%%%%%%
	\begin{dfn}
		For $x \in X$ and $r>0$, define open ball and closed ball with radius $r$ and center $x$ respectively  as: 
		$$ B_\phi (x, r)= \lbrace {y \in X : d_\phi(x, y) <r}\}~  \text{and}~ B_\phi [x, r]= \lbrace {y \in X : d_\phi(x, y) \leq r}\}. $$	 
	\end{dfn} 
	%%%%%%%%%
	\begin{ppn} Let $(X,d_\phi)$ be a $\phi$-metric space. Then for all $ r, s >0 $ and for all $ a\in X$,
		\begin{enumerate}[(i)]
			\item $r \leq s $ if and only if $ B_\phi(a,r)\subseteq B_\phi(a,s)$.
			\item $r \leq s $ if and only if $ B_\phi[a,r]\subseteq B_\phi(a,s)$. 
			\item $B_\phi(a,r)\subseteq B_\phi[a,r]$.
		\end{enumerate}
	\end{ppn}
	%%%
	%%%%%%%%%%%%%%%%%%%%%%%%%%%% TOPOLOGY PROOF %%%%%%%%%%%%%%%%%%%%%%%%%%%
	\begin{thm}
		Let $(X,d_\phi)$ be a $\phi$-metric space and define  $$ \tau_\phi=\{ G\subseteq X: \forall ~x\in G, \exists ~r>0 ~such~ that~ B_\phi(x,r)\subseteq G  \}.  $$
		Then $\tau_\phi$  is a topology on $X$.
		\begin{proof}
			Obviously $ \phi, X\in \tau_\phi$ and $\tau_\phi $ is closed under arbitrary union.  To check the closedness of $\tau_\phi $ under finite intersection, let us consider  $G_1, G_2\in \tau_\phi $. We need to show $G_1\cap G_2\in \tau_\phi $.\\ 
			Take any $x\in  G_1\cap G_2$. Then $x\in  G_1$ and $x\in  G_2$. So there exists  $ r_1, r_2 > 0 $ such that $B_\phi(x,r_1)\subseteq G_1 $ and $B_\phi(x,r_2)\subseteq G_2 $.\\
			Now if $r = \min\{ r_1, r_2\} $, then $B_\phi(x,r)\subseteq B_\phi(x,r_1)\subseteq G_1 $ and $B_\phi(x,r)\subseteq B_\phi(x,r_2)\subseteq G_2  $.\\
			Thus $B_\phi(x,r)\subseteq G_1\cap G_2 $. So, $G_1\cap G_2\in \tau_\phi $.	
		\end{proof}
	\end{thm}
	%%%%%%%%%%%%%%%
	\begin{dfn} Let $(X,d_\phi)$ be a $\phi$-metric space and $B\subseteq X$. Then
		\begin{enumerate}[(i)]
			\item  $B$  is said to be open set if $B\in \tau_\phi$.
			\item $B$ is said to be closed set if $X\setminus B$ is in $\tau_\phi$.
			\item $x\in X$ is called a limit point of $B$ if there exists $ r>0 $ such that $ { (B(x,r)\setminus \{x\}) }\cap {B }  $ contains infinitely many points of $B$. 
			\item The set of all limit points of $B$ is called the derived set of $B$, denoted by $B'$.
			\item A set that contains  the points of  $B$, as well as limit points of $B$, is called the closure of the set $B$, denoted by $\overline{B}$.
		\end{enumerate}
	\end{dfn}
	%%%%%%%%%%%%%%%
	The following two propositions in  $\phi$-metric space are obvious.
	%%%%%%%%%%%%%%%%%%%%
	\begin{ppn}Let $(X,d_\phi)$ be a $\phi$-metric space and $A\subseteq X$. Then
		\begin{enumerate}[(i)]
			\item $\overline{A}$ is a closed set.
			\item  $A$ is a closed set if and only if $A=\overline{A}$.
			\item  $x\notin \overline{A} $ if and only if $ d_\phi (x,a)>0,$ for all $ a \in A$.
		\end{enumerate}
	\end{ppn}
	%%%%%%%%%%%%%%%%%%%%%%%%%%%%%%%%%%%%%%%%%%%%
	\begin{ppn} In a $\phi$-metric space $(X,d_\phi)$,
		\begin{enumerate}[(i)]
			\item Arbitrary union of open set is open.
			\item  Arbitrary intersection of closed set is closed.
		\end{enumerate}
	\end{ppn}
	%%%%%%%%%%%
	\begin{rem}
		Arbitrary union (respectively intersection) of closed (respectively open) set is not closed (respectively open), which can be justified by examples of metric spaces.  
	\end{rem}
	%%%%%%%%%%%%%%%%%%%%%%%%%%%%
	Now we are interested to find a basis for $\tau_\phi $. In fact, we want to show that the set of all open balls form a basis. For the first step, we prove the next result.
	%%%%%%%%%%%%%%%%%%%%%%
	%%%%%%%%%%%%%%%%%%%%%%%%%%%%%
	%%%%%$$$$$$$$$$$$$$$$$$$$$%%% OPEN BALL OPEN SET %%%%%%%%%%%%%%%%%%%%%%%%%%%%%
	\begin{thm}
		In a $\phi$-metric space, every open ball is an open set.
		\begin{proof}
			For some $x\in X $ and $r>0 $, consider the open ball $ B_\phi(x,r)$ and choose  $y\in B_\phi(x,r) $. Then $d_\phi(x,y)=r'(say)<r $. \\
			We need to find some $s>0  $ such that $B_\phi(y,s)\subset B_\phi(x,r) $.\\
			Again for $ x\in X, ~y\in B_\phi(x,r) $, and $ a \in X, $ we have
			\begin{equation}\label{open}
				d_\phi(x,a)\leq  d_\phi(x,y)+d_\phi(y,a)+\phi(x,a,y)
			\end{equation}  
			Now for $ \dfrac{r-r'}{2}(>0) $, there exists  $ t>0  $ such that $ ~\phi(x,z,y)<  \dfrac{r-r'}{2} $ whenever $d_\phi(z,y)<t $, and $ z \in X $.\\
			Take  $ s =\min \{ \dfrac{r-r'}{2} ,t \}$. \\
			Let us choose $ z \in B_\phi(y,s) $. Then $ d_\phi(y,z) < s $, and hence $ \phi(x,z,y)<  \dfrac{r-r'}{2} $. \\
			Therefore from the inequality (\ref{open}) we have, 
			\begin{align*}
				d_\phi(x,z)\leq & d_\phi(x,y)+d_\phi(y,z)+\phi(x,z,y)\\
				< & r'+  s + \dfrac{r-r'}{2} \\ 
				\leq & r' + 2(\dfrac{r-r'}{2})=r
			\end{align*} 
			Thus, $d_\phi(x,z)<r  $ whenever $z\in B_\phi(y,s) $ where $ s =\min \{ \dfrac{r-r'}{2} ,t \} $.\\
			Hence $ B_\phi(y,s)\subset B_\phi(x,r) $, and this proves that $  B_\phi(x,r) $ is an open set.
		\end{proof}
	\end{thm}	
	%%%%
	%%
	%%%%%%%%%%
	%%%%%%%%%%%%%%%%%%%%%%%%%%%%%%%%%%%%%%%%%%%%%%%%%%% BASE FOR TOPOLOGY %%%%%%%%%%%%%%%%
	Consider the collection of open balls $ \beta=\{B_\phi(x,r): x\in X, r>0\} $. Now we will show that it generates a topology on $X$.
	\begin{thm}
		Let $(X, d_\phi)$ be a $\phi$-metric space. Then $\beta $ is a base for  $(X, \tau_\phi) $. 
		\begin{proof}
			Let $ x\in X$. Then for any $r>0,~x\in B_\phi(x,r) $. Next suppose for some $x $ and $y $ in $X $ and for $r_1,r_2 >0$, there is $z\in B_\phi(x,r_1)\cap B_\phi(y,r_2) $.\\
			Since $z\in B_\phi(x,r_1)$, so  there exists  $ s_1>0 $ such that $B_\phi(z,s_1)\subseteq B_\phi(x,r_1) $ and  $z\in B_\phi(y,r_2)$ implies  there exists  $ s_2>0 $ such that $B_\phi(z,s_2)\subseteq B_\phi(x,r_2).$\\
			Now take $s=\min\{ s_1, s_2 \} $. Then $z\in B_\phi(z,s)\subset B_\phi(x,r_1) \cap B_\phi(y,r_2)  $. 	
		\end{proof}
	\end{thm}
	%%%%%%%%%%%%%%%%%%%%%%%%%%%%%%%
	%%%%%%%%%%%%%%%%%%%%%%%%%%%%CLOSED BALL CLOSED SET %%%%%%%%%%%%%%%%%
	\begin{thm}
		In a $\phi$-metric space $(X,d_\phi)$, every closed ball is a closed set.
		\begin{proof}
			For any $x\in X $, and $r>0 $, consider the closed ball $ B_\phi[x,r]$.
			To prove that $ B_\phi[x,r]$ is closed, it is enough to show that $X \setminus B_\phi[x,r]=A$(say) is open.\\
			Choose $y\in A $. Then $ d_\phi(x,y) = r'(say)>r $.\\
			Now we need to find some $s>0  $ such that $B_\phi(y,s)\subset A $.\\
			For $(\dfrac{r'-r}{2})>0$, there exists  $ t>0  $ such that  $ \phi(x,y,z)< \dfrac{r'-r}{2} $ whenever $d_\phi(y,z)<t $, and $z\in X $. \\
			Let $s=\min \{ \dfrac{r'-r}{2}, t \} $.\\
			Choose   $ a \in B_\phi(y,s) $. Then $d_\phi(a,y) < s $, and so  $ \phi(x,y,a)< \dfrac{r'-r}{2} $. \\
			
			So for  $x\in X, ~y\in A $, and $ a\in B_\phi(y,s)$, we have
			\begin{align*}
				& d_\phi(x,y)\leq  d_\phi(x,a)+d_\phi(y,a)+\phi(x,y,a) \\
				or ~~ & d_\phi(x,a)  \geq d_\phi(x,y) -  d_\phi(y,a) - \phi(x,y,a)
			\end{align*}
			Hence,
			$$  d_\phi(x,a)  >   r' - s - (\dfrac{r'-r}{2}) \geq   r'- 2(\frac{r'-r}{2}) = r $$
			Therefore, $d(x,a)>r  $ whenever $ a \in B_\phi(y,s) $ where $ s=\min \{ \dfrac{r'-r}{2}, t \} $, which implies $A$ is an open set, and consequently $ B_\phi[x,r]$ is a closed set.
		\end{proof}
	\end{thm}
	%%%%%%   %%%%%%%%%%%%%%
	%%%%%%%%%%%%$$$$$$$$$$$$$
	%%%  regular   %%%%%%%%%%%%%%%%%%%
	\begin{thm}\label{regular}
		Every $\phi$-metric space  is  regular.
		\begin{proof}
			Let $ A $ be a closed in $(X,d_\phi)$ and $x\in X\setminus A$. So $d_\phi(x,a)>0,$ for all $ a\in A$.\\
			Let $3r=\inf\{ d_\phi(x,a):a\in A \} $. Consider the open ball  $ B(x,r)= V(say) $.\\  
			Now for any $a \in A$ and $r >0$, there exists  $ \beta>0 $ such that  $\phi(x,a,y)<r$ whenever $d_\phi(a,y)<\beta $, and $ y\in X $.\\
			Take $ \min \{\beta, r \}=r^*(say)$ and consider the open set $U=\underset{a\in A}{\cup} B_\phi (a,r^*) $. Then $ A \subset U $. \\ 
			We claim that  $U \cap  V=\phi$.  If possible suppose  there exists  $  c \in U\cap V  $. Then for any $ a\in A, ~d_\phi(a,c)< r^*  $, 
			and $ \phi(x,a,c) < r  $. \\
			Now for any $a\in A$, 
			$$  	d_\phi(x,a)\leq  ~ d_\phi(x,c)+ d_\phi(c,a)+\phi(a,x,c)	<  ~ r ^* + r + r \leq 3r. $$
			This is a contradiction to our assumption and hence $U$ and $V$ are two disjoint open sets in $X$ containing $A$, and $x$ respectively.
		\end{proof} 
	\end{thm}
	%%%%%%%%%%%%%%%%%%%%%%%%%%%%%%%%%%%%%%%%%%%%%%%%%%%%%%%%%
	%%%%%%%%%%%%%%%%%%%%%%%%%%%  NORMAl          %%%%%%%%%%%%%%%%%%%%
	\begin{thm}
		Every $\phi$-metric space  is normal.
		\begin{proof}
			Let $ A $ and $ B $ be two closed disjoint sets in $(X,d_\phi)$. Then for any $ a \in  A$, and $ b \in B,~ d_\phi(a,b) > 0 $.\\ 
			Let  $3r= \inf\{ d_\phi(a,b):a\in A,b\in B \} $. Consider the open set $V= \underset{b\in B}{\cup} B_\phi(b,r) $ containing $B$. \\
			For any $a\in A, ~b\in B$, and  $r>0  $, there exists  $ \delta>0 $ such that  $\phi(a,b,z) < r $ whenever $d_\phi(a,z)<\delta$ and $z\in X$.\\
			Let $ r^* = \min \{\delta, r \} $ and $ U=\underset{a\in A}{\cup} B_\phi (a,r^*) $. Then $ U $ is open and $A\subset U$. \\
			Next, we claim that  $U $ and $ V$ are disjoint. If not, then  there exists  $ c \in U\cap V  $. Then  for all $ a \in A $ and \\
			for all  $ b\in B, ~ d_\phi(a,c)< r^*, ~d_\phi(b,c)<r $ and $ \phi(a,b,c)<r $.\\
			So for $a\in A, ~b\in B$ and $c \in U\cap V $, 
			$$	d_\phi(a,b)\leq  ~ d_\phi(a,c)+d_\phi(c,b)+\phi(a,b,c)	<  ~ r^*+r+r\leq 3r. $$
			This contradicts our assumption and hence the theorem is proved.	
		\end{proof}
	\end{thm}
	%%%%%%%%%%%%%%%%%%%%%%%%%%%%%%%%%%%%%
	%%%%%%%%%%%%%%%%%%%%%%(metrizability)%%%%%%%%%%%%%%%%%%%%%%%%%%%%%%%%
	Now we  prove the Stone-type theorem in a $\phi$-metric space and use Bing metrization theorem to obtain a sufficient condition of metrizability.
	%%%%%%%%%%%
	\begin{thm}\label{thm-stone}  (Stone-type theorem) In a $ \phi- $metric space $ (X,d_\phi) $ every open cover of $X$ has an open refinement which is both $\sigma$-locally finite and $\sigma$-discrete.
		\begin{proof}
			Let $\{ \mathcal{U}_s:s \in S \}$ be an open cover of $X$. By the Zermelo theorem on well-ordering\cite{[10]}, we can take a well-ordering relation $<$ on $S$. Define the families
			$\mathcal{V}_i=\{\mathcal{V}_{s,i}:s \in S \}$ of subsets of $X$ by letting $\mathcal{V}_{s,i}= \underset{c\in \mathcal{C}}{\cup} B_\phi (c,\frac{1}{2^i}),$ where $\mathcal{C}$ is the set of all points $c \in X$ satisfying following conditions:
			\begin{enumerate}[(i)]
				\item $s$ is the smallest element of $S$ such that $c \in \mathcal{U}_s$.
				\item $c \notin \mathcal{V}_{t,j}$ for all $j<i$, and for all $t \in S$.
				\item $B_\phi(c, \frac{5}{2^i})\subset \mathcal{U}_s$.
			\end{enumerate}
			Obviously the sets $\mathcal{V}_{s,i}$ are open and by condition (iii), we have $V_{s,i}\subset \mathcal{U}_s$.\\
			For each $x \in X$, take the smallest $s\in S$ such that $x \in \mathcal{U}_s$, and a natural number $i$ such that $B_\phi(x,\frac{5}{2^i})\subset \mathcal{U} _s$. It implies that $x \in\mathcal{C}$ if and only if $x \notin \mathcal{V}_{t,j}$ for all $j<i$, and for all $t \in S$. Then $x \in \mathcal{V}_{s,i}$. Thus we have either $x \in \mathcal{V}_{t,j}$ for all $j<i$ and for all $t \in S$ or $x \in \mathcal{V}_{s,i}$. This proves that $\mathcal{V}=\underset{i\in \mathbb{N}}{\cup}\mathcal{V}_i$ is an open refinement of the cover  $\{ \mathcal{U}_s :s \in S \}$.\\
			Now  for every $i \in \mathbb{N}$, let $x_1\in \mathcal{V}_{s_1,i}$ and $x_2\in \mathcal{V}_{s_2,i}$ with $s_1\neq s_2$. Let us assume  $s_1< s_2$. By the definition of $V_{s_1,i}$ and $\mathcal{V}_{s_2,i} $, there exists  $ c_1,~ c_2\in X $ satisfying conditions (i), (ii), (iii) and $x_1\in B_\phi (c_1, \frac{1}{2^i}),~ x_2\in B_\phi (c_2, \frac{1}{2^i})$. Again we have $B_\phi(c_1, \frac{5}{2^i})\subset \mathcal{U}_{s_1}$ and $ c_2\notin \mathcal{U}_{s_1}$ and this implies  $d_\phi(c_1,c_2)\geq \frac{5}{2^i}$. But we have 
			\begin{equation*}
				d_\phi(c_1,c_2)\leq d_\phi(c_1,x_1) + d_\phi(x_1,x_2) + d_\phi(x_2,c_2) + \phi(c_1, c_2, x_1) +  \phi(x_1, c_2, x_2)
			\end{equation*}
			which implies 
			\begin{equation}\label{stone} 
				d_\phi(x_1,x_2) \geq \frac{5}{2^i} - d_\phi(c_1,x_1) - d_\phi(x_2,c_2) - \phi (c_1,c_2,x_1) - \phi (x_1,c_2,x_2).
			\end{equation}
			Again for $\frac{1}{2^{i+1}}(>0) $, there exists $ \beta_1, ~\beta_2 (>0)$ such that  $ \phi (c_1,c_2,x_1)< \frac{1}{2^{i+1}}$ whenever $ d_\phi(c_1, x_1)< \beta_1 $ and $ \phi (c_1,c_2,x_2)< \frac{1}{2^{i+1}}$ whenever $d_\phi(c_2, x_2)<\beta_2 $.\\
			Take $\min \{ \beta_1,~\beta_2,~ \frac{1}{2^i} \} = \beta(say)$. \\
			Then   $ d_\phi(c_1, x_1)< \beta, ~ d_\phi(c_2, x_2)< \beta $ and $  \phi (c_1,c_2,x_1)< \frac{1}{2^{i+1}}, ~ \phi (c_1,c_2,x_2)< \frac{1}{2^{i+1}} $.   Then (\ref{stone}) gives,
			$$  d_\phi (x_1, x_2)   > \frac{5}{2^i} -2\beta -2\times \frac{1}{2^{i+1}} 
			\geq \frac{5}{2^i} -2 \times \frac{1}{2^i}  -\frac{1}{2^{i}} = \frac{1}{2^{i-1}}. $$
			To prove that the families $V_i$ are $\sigma$-discrete, suppose  there exists $  x\in X$ such that 
			$ x_1, x_2 \in B_\phi (x, \frac{1}{2^{i+1}})$. Then we have  $  d_\phi(x, x_1) < \frac{1}{2^{i+1}}, ~  d_\phi(x, x_2) < \frac{1}{2^{i+1}}  $ and 
			\begin{equation}\label{stone2}
				\frac{1}{2^{i-1}}   < d_\phi (x_1, x_2)\leq  d_\phi (x_1, x) +  d_\phi (x, x_2) +\phi (x_1,x_2,x).
			\end{equation} 
			Now for $\frac{1}{2^{i+1}}>0 $, there exists  $ \beta'>0$ such that  $ \phi (x_1,x_2,x)< \frac{1}{2^{i+1}}$ whenever $ d_\phi(x_2, x) < \beta' $. \\
			If $ \delta=\min\{ \beta', \frac{1}{2^{i+1}} \}$ then  $ d_\phi(x_2, x) < \delta $ and  $ \phi (x_1,x_2,x)< \frac{1}{2^{i+1}}  $. The inequality (\ref{stone2}) gives, 
			$$ \frac{1}{2^{i-1}}   < d_\phi (x_1, x_2)<  \frac{1}{2^{i+1}} + \delta + \frac{1}{2^{i+1}} 
			\leq  2\times \frac{1}{2^{i+1}} + \frac{1}{2^{i+1}} 
			<   \frac{1}{2^i} + \frac{1}{2^{i}} = \frac{1}{2^{i-1}}. $$
			This is a contradiction and hence it proves that each ball of radius $\frac{1}{2^{i+1}}$ meets at most one element of $\mathcal{V}_i$ that is $\mathcal{V}=\cup _{i\in \mathbb{N}} \mathcal{V}_i$ is $\sigma$-discrete.\\
			Let $i\in \mathbb{N}$, then for all $t\in S,~i\geq j+k  $ and $c\in\mathcal{C}$ implies $c\notin \mathcal{V}_{t,j}$. Now if $B_\phi (x,\frac{1}{2^{k-1}})\subset \mathcal{V}_{t,j}$, then $c \notin B_\phi (x,\frac{1}{2^{k-1}})$ and $d_\phi(x,c)\geq \frac{1}{2^{k-1}}$. Again $j+k\geq k+1$ and $i\geq k+1$ implies $\frac{1}{2^{j+k}}\leq \frac{1}{2^{k+1}}$ and $\frac{1}{2^i}\leq \frac{1}{2^{k+1}}$. \\
			Next suppose there exists  $ y\in B_\phi (x,\frac{1}{2^{j+k}})\cap B_\phi (c,\frac{1}{2^i})$. Then 
			\begin{equation}\label{stone3}
				d_\phi(x,c) \leq d_\phi(x,y) + d_\phi(y,c) + \phi(x,c,y).
			\end{equation}
			For $\frac{1}{2^k}>0 $, there exists $  \alpha>0$ such that  $ \phi (x,c,y)< \frac{1}{2^k}$ whenever $d_\phi(x,y)<\alpha $.\\
			Let $ \gamma=\min\{ \alpha, \frac{1}{2^{j+k}} \} $. Then $ d_\phi(x,y) < \gamma $ and $ \phi (x,c,y) < \frac{1}{2^{k}} $.  Therefore from (\ref{stone3}) we obtain,
			$$
			\frac{1}{2^{k-1}}  \leq d_\phi(x,c)   < \gamma  + \frac{1}{2^i} + \frac{1}{2^k} 
			\leq \frac{1}{2^{j+k}} + \frac{1}{2^{k+1}} + \frac{1}{2^k}
			\leq \frac{1}{2^{k+1}} + \frac{1}{2^{k+1}} + \frac{1}{2^{k}}
			= \frac{1}{2^{k-1}}. $$
			which concludes  $B_\phi (x,\frac{1}{2^{j+k}})\cap B_\phi (c,\frac{1}{2^i})=\phi$, and this implies $B_\phi (x,\frac{1}{2^{j+k}})\cap \mathcal{V}_{s,i}=\phi$ for $i \geq j+k$ and $s\in S$ with $B_\phi (x,\frac{1}{2^{k-1}})\subset \mathcal{V}_{t,j}$. \\
			Again for each $x \in X$, since $\mathcal{V}$ is a refinement of $\{ \mathcal{U}_s:s\in S \} $, there exists $ l, j $, and $t$ such that $B_\phi (x,\frac{1}{2^l})\subset \mathcal{V}_{t,j}$ and thus   there exists $ k, j$, and $t$ such that $B_\phi (x,\frac{1}{2^{k-1}})\subset \mathcal{V}_{t,j}$. Then the ball $B_\phi (x,\frac{1}{2^{j+k}})$ meets at most $(j+k-1)$ members of $\mathcal{V}$. This proves that $\mathcal{V}_i$ is locally finite that is $\mathcal{V}$ is $\sigma$-locally finite.
		\end{proof}
	\end{thm} 
	%%%%%%%%%%%%%%%%
	\begin{cor}\label{cor-sigma discrete}
		Let $ (X,d_\phi) $ be a $ \phi$-metric space. Then $X$ has $\sigma$-discrete base.
		\begin{proof}
			For every $i \in \mathbb{N}$, let $\mathcal{A}_i=\{ B_\phi(x,\frac{1}{i}) :x \in X \}$. Then $\mathcal{A}_i$ is an open cover of $X$. By Theorem \ref{thm-stone}, there exists an open $\sigma$-discrete refinement $\mathcal{B}_i$ of $\mathcal{A}_i$. Put $\mathcal{B}=\underset {i\in \mathbb{N}}{\cup}\mathcal{B}_i$. Then $\mathcal{B}$ is a $\sigma$-discrete base of $X$.
		\end{proof}
	\end{cor}
	%%%%%%%%%%%%%%%%%
	\begin{cor}
		Every  $ \phi$-metric space is metrizable.
		\begin{proof}
			From Theorem \ref{regular} and Corollary \ref{cor-sigma discrete} it follows that $X$ is regular space with $\sigma-$discrete base. Then from Bing metrization theorem(Theorem \ref{thm-bing}), $X$ is metrizable.
		\end{proof}
	\end{cor}
	%%%%%%%%%%%%%%%%%%%%%%%%%%%%
	%%%%%%%%%%%%%%%%%%%%%%%%%%%%%%%%%Hausdorff%%%%%%%%%%%%%%%
	\begin{cor}
		Let $(X,d_\phi)$ be a $\phi$-metric space and $\tau_\phi$ is a topology on $X$. Then $\tau_\phi$ is a Hausdorff topology on $X$. 
	\end{cor}
	\begin{proof}
		Since $(X,\tau_\phi)$ is a regular topological space, so it is Hausdorff.
	\end{proof}
	%%%%%%%%%%%%%%
	\begin{rem}
		Till now we have shown that $\phi$-metric can be induced from a b-metric, S-metric, etc. From this, we deduce that open balls of S-metric(or b-metric) are the same as the open balls of induced $\phi$-metric. Therefore  topology generated by the open balls of S-metric(or b-metric) is identical to the topology generated by the open balls of respective induced $\phi$-metric. Thus S-metric and  b-metric spaces are $\phi$-metrizable.
	\end{rem}
	%%%%%%%%%%%%%%%%
	%%%%%%%%%%%%%%%%%%%%%%
	%%%%%%%%%%%%%%%%
	%%%%%%%%%%%%%%%%%%%%%%%%%%
	%%%%%   convergence of sequence%%%%%%%%%
	%%%%%
	Next, we  discuss the convergence of a sequence in $\phi$-metric space including its basic properties.
	%%%%%%%%%%%%%%%%%%%
	%%% CONVERGENCE DEFN %%%%%%%%%%%%%%%%%%%%%%%%%%%%%%%%%%
	\begin{dfn}
		A sequence $\lbrace x_n \rbrace\subseteq X$ is said to converge to $x$ if for any $\epsilon>0,$ there exists a positive integer  $ N $  such that 
		$$ d_\phi(x_n,x)< \epsilon,~ \text{for all } ~ n \geq N ~~~\text{ that is} ~ d(x_n,x)\rightarrow 0~ \text{ as}~  n \rightarrow \infty. $$ 
		We denote this by $ \underset{n\rightarrow \infty}{lim}~ x_n = x$.	
	\end{dfn}
	%%%%%%%%%%%%%%%%%%% UNIQUE LIMIT %%%%%%%%%%%%%%%%%%%%%%%%%%%%%%%%%%%%%%%
	\begin{ppn}
		In a $\phi$-metric space $(X,d_\phi)$, every convergent sequence has unique limit.
		\begin{proof}
			Since $(X,\tau_\phi)$ is a Hausdorff topological space, so the conclusion is obvious.
		\end{proof}
	\end{ppn}
	%%%%%%%%%%%%%%%%%%%%%%%%%%%
	\begin{ppn}
		Let $(X,d_\phi)$ be a $\phi$-metric space, and $d$ be the metric on $X$ whose topology is identical to the $\phi$-metric topology. Then for any sequence  $\{x_n\}\subseteq X $, and $x\in X$,
		$$ \underset{n\rightarrow\infty}{lim} d_\phi(x_n,x)=0 ~ \text{if and only if } ~ \underset{n\rightarrow\infty}{lim}d(x_n,x)=0. $$
		\begin{proof}
			First assume that $~\underset{n\rightarrow\infty}{lim} d_\phi(x_n,x)=0$. Then, for all $ \epsilon >0 $, there exists $ N\in \mathbb{N} $ such that
			$$  d_\phi(x_n,x)<\epsilon, ~ \text{for all} ~ n\geq N ~~ that ~ is ~~  x_n\in B_\phi(x,\epsilon),~~\text{for all} ~n\geq N  $$
			Hence, there exists $ \delta(\epsilon)>0 ~\text{such that}~ x_n\in B_\phi(x,\epsilon)\subset B(x,\delta), ~ \text{for all} ~ n\geq N $. \\
			Therefore, $ \underset{n\rightarrow\infty}{lim}d(x_n,x)=0. $ \\
			If conversely assume  $~\underset{n\rightarrow\infty} {lim} d(x_n,x)= 0 $, then for all  $ \epsilon >0,$ there exists $  N\in \mathbb{N} $  such that
			\begin{align*}
				& d(x_n,x) < \epsilon,~~\text{for all} ~ n\geq N  \\
				or ~~ & x_n\in B(x,\epsilon),~ \text{for all} ~  n\geq N  
			\end{align*}
			Hence, there exists $ \delta(\epsilon)>0 ~\text{such that}~ x_n\in B(x,\epsilon)\subset B_\phi(x,\delta) ,~~\text{for all} ~ n\geq N  $. \\
			Thus,	$  \underset{n\rightarrow\infty}{lim}d_\phi(x_n,x)=0. $ \\
			Hence the proof is complete.
		\end{proof}	 	
	\end{ppn}
	%%%%%% 
	%%%%%%%%%%%%
	%%%%%%%%
	\begin{ppn} 
		Let $(X,d_\phi)$ be a $\phi$-metric space, and $\{x_n\}$ and $\{y_n\} $ be two sequences in  $X $ converging to  $ x $ and  $ y $ respectively. Then the sequence $\{ d_\phi(x_n,y_n)\} $ converges to $ d_\phi(x,y)$. 
		\begin{proof}
			Let $\epsilon>0$. We have,
			\begin{align*}
				d_\phi(x,y) & \leq d_\phi(x,x_n)+ d_\phi(x_n,y)+\phi(x,y,x_n)\\
				&\leq d_\phi(x,x_n)+ d_\phi(x_n,y_n)+ d_\phi(y_n,y)+\phi(x_n,y,y_n) +\phi (x,y,x_n)
			\end{align*} 
			and
			\begin{align*}
				d_\phi(x_n,y_n)   \leq &  d_\phi(x_n,x)+ d_\phi(x,y_n) + \phi( x_n, y_n, x)\\
				\leq &  d_\phi(x_n,x)+ d_\phi(x,y)+ d_\phi(y,y_n)+  \phi(x,y_n,y) + \phi(x_n,y_n,x)
			\end{align*}
			Now for  $\frac{\epsilon}{4}(>0) $, there exists $ \beta_1, \beta_2 >0$ such that $\phi(z,y,y_n)<\frac{\epsilon}{4}, ~ \phi( z,y_n,y) < \frac{\epsilon}{4} $ whenever $ d_\phi(y_n,y)<\beta_1,\\
			~z\in X $ and  $\phi (x,w,x_n)<\frac{\epsilon}{4}, ~ \phi( x_n, w, x) < \frac{\epsilon}{4}$ whenever $ d_\phi(x_n,x)<\beta_2, ~w\in X $.\\
			Let $\delta=\min \{\beta_1, \beta_2, \frac{\epsilon}{4}\}$. \\
			Then for  $\delta>0$, there exists $ N_1, N_2\in \mathbb{N} $ such that $ d_\phi(x_n,x)<\delta $, for all $ n\geq N_1 $ and $ d_\phi(y_n,y)<\delta $, for all $ n\geq N_2 $. \\ 
			Take $\max\{N_1,N_2\}=N(say)$. \\
			Then  for all $ n \geq N $, 
			$ d_\phi(x_n,x)<\delta $ and $ d_\phi(y_n,y)<\delta ~$ implies  $\phi (x,y,x_n)< \frac{\epsilon}{4}, ~\phi(x_n,y,y_n)<\frac{\epsilon}{4}, \\
			~\phi( x_n, y_n, x) <\frac{\epsilon}{4} $ and $ \phi(x,y_n, y)< \frac{\epsilon}{4} $.\\ 
			Thus for all $ n\geq N$ we have,
			$$	d_\phi(x,y)  <  \delta + d_\phi (x_n,y_n) +\delta  +\frac{\epsilon}{4} + \frac{\epsilon}{4} 
			\leq 4\cdot \frac{\epsilon}{4} + d_\phi (x_n,y_n)  
			= \epsilon +d_\phi(x_n,y_n) $$
			and 
			$$	d_\phi(x_n,y_n) <  \delta + d_\phi(x,y) + \delta +\frac{\epsilon}{4} + \frac{\epsilon}{4} 
			\leq 4\cdot \frac{\epsilon}{4} + d_\phi (x,y)  
			= \epsilon +d_\phi(x,y) $$
			Since $\epsilon > 0$ is arbitrary, by taking limit as $n\rightarrow \infty$ on both side, we obtain $d_\phi(x,y) \leq \underset{n\rightarrow\infty }{lim} d_\phi(x_n,y_n)  $, and $\underset{n\rightarrow\infty}{lim} d_\phi(x_n,y_n)  \leq d_\phi(x,y) $, which implies $\underset{n\rightarrow\infty}{lim} d_\phi(x_n,y_n)=d(x,y)$.	
		\end{proof} 
	\end{ppn}
	%%%%%%%%%%%% CAUCHY %%%%%%%%%%%%%%%%%%%%%%%%%%%%%%%%%%%
	\begin{dfn}
		In a $\phi$-metric space $(X,d_\phi)$, a sequence $\lbrace x_n \rbrace\subseteq X$ is said to be Cauchy if for any  $\epsilon>0,$ there exists a positive integer $ N $ such that 
		$$ d_\phi(x_n,x_m)< \epsilon, ~ \text{for all} ~ m,n\geq N ~~~ \text{ that is} ~ d_\phi(x_n,x_m)\rightarrow 0~ \text{ as} ~  m,n\rightarrow\infty. $$ 
	\end{dfn}
	%%%%%%%%%%%%
	%%%%%%%%%%%%% CONV IMPLY CAUCHY %%%%%%%%%%%%%%%%%%%%%%%%%%%%
	\begin{ppn} \label{ conv cauchy}
		In a $\phi$-metric space $(X,d_\phi)$, every convergent sequence is Cauchy.
		\begin{proof}
			Let $\epsilon >0$ and $\{x_n\}\subseteq X$ converges to $x$.\\
			Now, 
			$$ d_\phi(x_m,x_n)\leq d_\phi(x_m,x)+d_\phi(x,x_n)+\phi(x_m,x_n,x), ~ \text{for all} ~ m,n \in  \mathbb{N} $$ 
			For  $\frac{\epsilon}{3} $, there exists $ \beta>0$ such that $\phi(x_m,x_n,x)<\frac{\epsilon}{3}$ whenever  $d_\phi(x_m,x)<\beta $. \\
			Let $ \delta = \min\{ \beta, \frac{\epsilon}{3}\}$. \\
			Again for $ \delta>0,$ there exists a natural number $ N $ such that $d_\phi(x_n,x)< \delta $, for all $ n\geq N $.\\ 
			So for all $ m, n\geq N, ~ d_\phi(x_m,x)<\delta,  d_\phi(x_n,x)< \delta $ and  $\phi(x_m,x_n,x)<\frac{\epsilon}{3}$. 
			Thus for all $ m, n\geq N $,
			$$  d_\phi(x_m,x_n)   < \delta + \delta +  \frac{\epsilon}{3} 
			\leq 2\cdot \frac{\epsilon}{3} + \frac{\epsilon}{3}= \epsilon. $$ 
			Hence $\{x_n\}$ is a Cauchy sequence in $X$. 
		\end{proof}
	\end{ppn}	
	%%%%%%%%%%%%%%%%%%%%%%%%       %%%%%%%%%%%%%%%%%%%%%%%%%%
	As in metric space, we can define the boundedness of a set in a  $\phi$-metric space.
	%%%%%%%%%%%%%%%%%%%
	\begin{dfn}
		Let $(X,d_\phi)$ be a $\phi$-metric space.  $A \subseteq X $ is said to be bounded if there exists a non-negative real number $ K  $ such that $ d_\phi(x,y) < K $, for all $ x,y \in A  $.
	\end{dfn} 
	%%%%%%%%%%%%%%%%%%%%%%%
	\begin{ppn} \label{conv bdd}
		Every convergent sequence a $\phi$-metric space $(X,d_\phi)$, is bounded.
		\begin{proof}
			Let $\epsilon >0 $ and $ \{x_n\} $ be a sequence in $X$ converging to $x \in X $. \\
			Again we have, $$ d_\phi(x_n,x_m) \leq d_\phi(x_n,x) + d_\phi(x_m,x) + \phi(x_n,x_m,x), ~ \text{for all} ~ m,n \in \mathbb{N} $$
			Choose $\epsilon = 1$. \\
			For $ \epsilon = 1 $, there  exists $ \beta >0$  such that $  \phi(x_n,x_m,x) < 1 $ whenever $ d_\phi(x_m,x) < \beta $.  \\
			Let  $\delta= \min\{ \beta, 1 \}$. \\
			Then  for $\delta > 0 $, there exists $ N \in \mathbb{N}$ such that  $ d_\phi(x_m,x)< \delta $, for all $ m \geq N $. \\
			So for $ m, n \geq N, ~ d_\phi(x_n,x)< \delta, ~ d_\phi(x_m,x)< \delta $ and  $ \phi(x_n,x_m,x) < 1 $.
			Then  for all $ m, n \geq N $,  
			$$	d_\phi(x_n,x_m) < \delta + \delta +1 \leq 3. $$
			Now suppose $ M= \max \{ d_\phi(x_r,x_s):1\leq r,s < N \} $. Then $ d_\phi(x_n,x_m) \leq M $,  for all $ m, n < N $. \\
			If $K=max\{ M,3\}$ then $ d_\phi(x_n,x_m) \leq K $, for all $ m,n \in  \mathbb{N}  $.\\
			This completes the proof. 
		\end{proof}
	\end{ppn}
	%%%%%%%%%%%%%%%%%%%%%%%%%%%%%%
	\begin{rem}
		The converse of Proposition \ref{ conv cauchy} and Proposition \ref{conv bdd} are not true in general,  since those statements do not hold for metric spaces.   
	\end{rem}
	%%%%%%%%%%%%%
	\begin{rem}
		As open balls of an S-metric (or b-metric) space are the same as the open balls of induced $\phi$-metric space, so the convergence of a sequence 
		is also identical in both cases.
	\end{rem}
	%%%%
	%%%
	%%
	%%%%%%%%%%%%%
	%%%%%%%%%%Compactness -- completeness %%%%%%%%5
	%%%%%%%%%%%%%%%%%%%
	\section{Some basic properties of $\phi$-metric spaces}
	In Example \ref{ex2} we have shown that the product of a finite number of $\phi$-metrics on a non-empty set is also a $\phi$-metric on that set. 
	Now if we consider the set of all $\phi$-metrics on a non-empty set $X$ with the binary operation multiplication$(\cdot)$ then we obtain an 
	algebraic structure, say  $(\mathscr{D},\cdot)$. \\
	Let us discuss its structure in detail. \\ 
	Already we have shown $\mathscr{D}$ is closed under $`\cdot$' and obviously  $`\cdot$' is both commutative and associative on $\mathscr{D}$. 
	So $(\mathscr{D},\cdot)$ forms a commutative semigroup. The discrete metric is the only idempotent element and there is no nilpotent element in 
	$(\mathscr{D},\cdot) $.\\
	Moreover $(\mathscr{D},\cdot)$ is a commutative monoid as the discrete metric on $X$ acts as  the identity element.\\
	~\\
	%%%%%%%%%%%%%%%%%%%%%%%%
	Next, we are interested to study the topology generated by a finite number of $\phi$-metrics. It is enough to study for the product of only two 
	$\phi$-metrics.
	%%%%%%%%%%%%%%%%%%%%%%%%%%%%%%%%%%%%%%%%%%%%%%%%
	\begin{ppn}
		Consider the $\phi$-metrics $D_\phi,~d_{\phi_1},~d_{\phi_2}$ where $ D_\phi=  d_{\phi_1}\cdot d_{\phi_2} $ and
		$\tau_\phi,~\tau_{\phi_1},~\tau_{\phi_2}$ be the topologies induced by the open balls of $D_\phi,~d_{\phi_1}$ and $d_{\phi_2}$ respectively. 
		Then $\tau_\phi=\tau_{\phi_1}\cap\tau_{\phi_2}$.
		\begin{proof}
			Let us denote the open balls of $(X,d_{\phi_i})$ by $B_{\phi_i}(x,r),~i=1,2$ and that of $(X,D_\phi)$ by $B_\phi(x,r)$, for some $x\in X,~r>0$. \\
			Choose $x\in X$, and $r>0$ such that $y\in B_{\phi_1}(x,r)\cap B_{\phi_2}(x,r)$. \\
			Then, $ d_{\phi_i}(x,y)<r,~ \text{for} ~ i=1,2 $. Therefore,
			$$  D_\phi(x,y) < r^2 ~~  that ~ is ~ ~ y\in B_\phi(x,r^2). $$
			This implies $ \tau_{\phi_1}\cap\tau_{\phi_2}\subseteq \tau_\phi. $ \\ 
			Again if $\{ x_n \}$ converges to $x$ in $(X,d_{\phi_1})$ then $\{ x_n \}$ converges to $x$ in $(X,D_\phi)$. So, for all $ \epsilon >0 $, there exists two positive integers $ N_1, N_2  $ such that 
			\begin{align*}
				&  d_{\phi_1}(x_n,x)<\epsilon,~ \text{for all} ~ n\geq N_1 ~ 
				\text{implies}~ D_\phi(x_n,x)<\epsilon,~ \text{for all} ~ n\geq N_2\\
				or ~~ & \text{for all} ~ n\geq N,~ x_n\in B_{\phi_1}(x,\epsilon) ~\text{implies}~  x_n\in B_{\phi}(x,\epsilon) ~~\text{where}~N= \max\{N_1,N_2\}
			\end{align*}
			This implies $ \tau_{\phi}\subseteq \tau_{\phi_1} $. \\
			Similarly,  $ \tau_{\phi}\subseteq \tau_{\phi_2} $.
			Therefore,  $ \tau_\phi\subseteq \tau_{\phi_1}\cap\tau_{\phi_2} $\\
			This completes the proof.	
		\end{proof}  
	\end{ppn}
	%%%%%%%%%%%%%%%%%%%%%%%%%%
	In the previous section, we have defined a bounded set in which the  distance between two elements is finite. This leads us to define the 
	diameter of a set and encourages us to check the relation between the  diameter of a set and its closure. 
	%%%%%%%%%%%%%%%%%%%%%
	\begin{dfn} Let $(X,d_\phi)$ be a $\phi$-metric space. Diameter of a set $F\subseteq X $ denoted by $\delta(F)$ and defined by  
		$\delta(F)=\underset{x,y\in F }{Sup}~ d_\phi(x,y)$.\\
		Therefore $F$ is said to be bounded, if $\delta(F)<\infty$, otherwise unbounded.
	\end{dfn}
	%%%%%%%%%%%%%%%%%%%%%%%%%%
	\begin{thm}
		For a subset $A$ of a $\phi$-metric space $(X,d_\phi),~  \delta(\overline{A})=\delta(A)$ where $\overline{A}$ denotes closure of $A$.
		\begin{proof}
			Since $A\subseteq \overline {A}$, so  
			\begin{equation} \label{ A A bar1}
				\delta(A)\leq \delta( \overline {A})
			\end{equation}
			Next choose any $\epsilon>0$, and $x,y\in \overline {A}$.  \\
			Now for $x, y\in \overline {A}$, and $ a, b \in X $, 
			\begin{align*}
				d_\phi(x,y) & \leq d_\phi(x,a)+ d_\phi(a,y)+\phi(x,y,a)\\
				&\leq d_\phi(x,a)+ d_\phi(a,b)+ d_\phi(b,y)+\phi(a,y,b) +\phi (x,y,a)
			\end{align*}
			Again for  $\frac{\epsilon}{4} $, there exists $ \beta_1, ~ \beta_2 >0$ such that $\phi(a,y,b)<\frac{\epsilon}{4}$ whenever $ d_\phi(b,y)<\beta_1$,
			and  $\phi (x,y,a)<\frac{\epsilon}{4}$ whenever $ d_\phi(a,x)<\beta_2$.\\
			Let $\gamma= \min \{\beta_1, \beta_2, \frac{\epsilon}{4}\}$.\\
			Since $x,y\in \overline {A}$, so there exists $ x_1\in A\cap B_\phi(x,\gamma)$, and $y_1\in A\cap B_\phi(y,\gamma) $.\\
			Therefore $ d_\phi(x,x_1) < \gamma, ~ d_\phi(y_1,y) < \gamma  $, and hence  $ \phi(x,y,x_1)< \frac{\epsilon}{4}, 
			~ \phi(x_1,y,y_1) < \frac{\epsilon}{4} $.\\
			Hence for $x,y\in \overline {A}$, and $ x_1\in A\cap B_\phi(x,\gamma), ~ y_1\in A\cap B_\phi(y,\gamma) $, 
			\begin{align*}
				d_\phi(x,y) & \leq d_\phi(x, x_1)+ d_\phi(x_1, y_1)+ d_\phi(y_1,y)+\phi(x_1,y,y_1) +\phi (x,y,x_1) \\
				& < \gamma+ d_\phi(x_1,y_1)+ \gamma + \frac{\epsilon}{4} + \frac{\epsilon}{4}\\
				& \leq 4\cdot\frac{\epsilon}{4}  +d_\phi(x_1,y_1)  
				\leq \epsilon + \underset{a,b\in A }{Sup}~ d_\phi(a,b)
			\end{align*}
			Since $\epsilon>0  $ is arbitrary, thus finally  we obtain
			$$ \underset{a,b\in \overline {A} }{Sup}~ d_\phi(a,b)\leq \underset{a,b\in A }{Sup}~ d_\phi(a,b)$$ 
			that is 
			\begin{equation}\label{ A A bar2}
				\delta( \overline {A})\leq \delta(A)
			\end{equation} 
			The relations (\ref{ A A bar1}) and (\ref{ A A bar2}) together gives, $ ~ \delta(\overline {A})= \delta(A)$.
		\end{proof}
	\end{thm}
	%%%%%%%%%%%%%%%%%%%%
	We now discuss  compactness, the most useful notion of a topological space including completeness.
	%%%%%%%%%%%%%%%%%%%%%%%%%%%%%%%%%
	\begin{thm}
		A $\phi$-metric space $(X,d_\phi)$ is compact if and only if it is sequentially compact.
		\begin{proof}
			Since $(X,d_\phi)$ is metrizable so there exist a metric on $X$, say $d$ whose topology is identical with the $\phi$-metric topology. \\
			Then $(X,d_\phi)$ is compact if and only if $ (X,d)$ is compact
			if and only if  $(X,d)$ is  sequentially compact 
			if and only if  $(X,d_\phi)$ is sequentially compact.
		\end{proof}
	\end{thm}
	%%%%%%%%%%%%%%%%%%%%%%%%%
	\begin{thm}
		Every compact $\phi$-metric space $(X,d_\phi)$ is closed and bounded.
	\end{thm}
	\begin{proof}
		If possible, suppose $X$ is not closed. So there exists a sequence of points $\{x_n\}$ such that $x_n\in X$ converges to a point $x\notin X$.\\
		Since $ X $ is compact, $\{x_n\}$ has  a subsequence which converges to a point in $X$. But subsequence must converge to $x$ which does 
		not belong to $X$. This contradicts the compactness of $X$. Hence $X$ is closed.\\
		Next, we prove that $X$ is bounded. \\
		If possible suppose that $X$ is unbounded, and choose $x_0\in X$, any fixed element. \\
		Since $X$ is unbounded, there exist $ x_1\in A$ such that $ d_\phi(x_1,x_0)>1 $. Similarly, there exists $ x_2\in X$ such that $ d_\phi(x_2,x_0)>2 $.
		Continuing in this way, there exists $ x_n\in X$ such that 
		$ d_\phi(x_n,x_0)>n $, for all $ n \in \mathbb{N} $.\\ 
		Since $X$ is compact, so there exists a subsequence  $\{ x_{n_k}\}$ of $\{ x_n\}$ such that
		$\underset{n\rightarrow\infty }{lim}  x_{n_k}=x\in X$. But we have, $d_\phi(x_{n_k},x_0)>n_k$. Again,
		\begin{equation}\label{compact is bdd}
			d_\phi(x_{n_k},x_0)\leq d_\phi(x_{n_k},x) + d_\phi(x,x_0) + \phi(x_{n_k},x_0,x)
		\end{equation}
		Let $ \epsilon >0$. \\
		Now for $ \epsilon >0,$ there exists $ \delta >0 $ such that $ \phi(x_l,x_0,x) < \frac{\epsilon}{2} $ whenever $ d_\phi(x_l,x) < \delta $.\\
		Let $ \beta=\min \{ \frac{\epsilon}{2}, \delta \} $. \\
		Since $ x_{n_k} \rightarrow x $ as $ k \rightarrow \infty$, so there exists $  N\in \mathbb{N} $ such that 
		$ d_\phi(x_{n_k},x) < \beta $, for all $ k \geq N $.\\
		Therefore for all $ k \geq N,~ d_\phi(x_{n_k},x) < \beta $, and $ \phi(x_{n_k},x_0,x) < \frac{\epsilon}{2} $.\\ 
		Thus the relation (\ref{compact is bdd}) gives, for all $ k \geq N $,
		$$	n_k < d_\phi(x_{n_k},x_0)  < \beta + d_\phi(x_0,x) + \frac{\epsilon}{2} 
		\leq  \frac{\epsilon}{2}  + d_\phi(x_0,x) +  \frac{\epsilon}{2} 
		= \epsilon + d_\phi(x_0,x)  $$
		Taking limit as $ k \rightarrow \infty $  on both side of the above inequality, we obtain $ \infty \leq  d_\phi(x, x_0) $. This  contradicts 
		that $d_\phi$ is a real valued function.\\
		Hence $X$ is bounded.	
	\end{proof}
	%%%%%%%%%%%%%%%%%%%%%%%%%
	%%%%%%%%%%%%%%%%%%%%%%% completeness %%%%%%%%%%%%%%%%%%%%%%%%
	\begin{dfn}
		A $\phi$-metric space $(X,d_\phi)$ is said to be complete if every Cauchy sequence in $X$ converges to some point in $X$.
	\end{dfn}
	%%%%%%%%%%%%
	\begin{thm}
		Every compact $\phi$-metric space $(X,d_\phi)$ is complete. 
	\end{thm}
	\begin{proof}
		Let $\epsilon >0$, and $\{ x_n\}$ be a Cauchy sequence in the compact $\phi$-metric space $(X,d_\phi)$.  \\
		So there exists a subsequence  $\{ x_{k_n}\}$ of $\{ x_n\}$ such that
		$ \underset {n\rightarrow\infty }{lim}  x_{k_n}= x\in X$.\\
		Now we have,  
		\begin{align*}
			d_\phi(x_n,x)\leq & ~ d_\phi(x_n,x_l) + d_\phi(x_l,x)  + \phi(x_n, x, x_l )\\
			\leq & ~ d_\phi(x_n,x_l) + d_\phi(x_{k_m},x_l) + d_\phi(x_{k_m},x) +  \phi(x_l, x, x_{k_m} ) + \phi(x_n, x, x_l )
		\end{align*}   
		For  $ \frac{\epsilon}{5} >0 $ there exists $ \delta_1,~\delta_2 >0$  such that $ \phi(x_n, x, x_l ) <\frac{\epsilon}{5}$ whenever $ d_\phi(x_n,x_l) <\delta_1 $, and  $ \phi(x_l, x, x_{k_m} ) <\frac{\epsilon}{5}$ whenever $ d_\phi(x_{k_m},x) <\delta_2 $.\\
		Let $\delta= \min \{\delta_1, \delta_2, \frac{\epsilon}{5}\}$.\\
		Since $\{ x_n\} $ is a Cauchy sequence,  so for  $\delta>0$, there exists a positive integer $n_0$ such that 
		$$ d_\phi(x_n,x_m)< \delta,~~\forall ~n,m\geq n_0 $$
		In particular,
		\begin{equation}\label{cauchy n0}
			d_\phi(x_n,x_{n_0})< \delta, ~~\text{for all} ~ n \geq n_0 
		\end{equation}
		Again $ x_{k_n} \rightarrow x  $ as $ n\rightarrow \infty $ implies $ \exists m \in \mathbb{ N}$ such that 
		\begin{equation}
			d_\phi(x_{k_m},x)< \delta,~~\text{for all} ~ m \geq n_0 
		\end{equation}
		Since $ k_m\geq m\geq n_0 $, so from (\ref{cauchy n0}), 
		\begin{equation}
			d_\phi(x_{k_m},x_{n_0})< \delta
		\end{equation}
		Therefore for $ n\geq n_0$, $ d_\phi(x_n,x_{n_0})< \delta, ~ d_\phi(x_{k_n},x)< \delta $, and $  \phi(x_n, x, x_{n_0} ) <  \frac{\epsilon}{5}, ~  \phi(x_{n_0}, x, x_{k_m} ) < \frac{\epsilon}{5}  $.\\
		So $\forall n\geq n_0$,
		\begin{align*}
			d_\phi(x_n,x) & \leq ~d_\phi(x_n,x_{n_0}) + d_\phi(x_{k_n},x_{n_0}) + d_\phi(x_{k_n},x) +  \phi(x_{n_0}, x, x_{k_n} ) + \phi(x_n, x, x_{n_0} )\\
			& < ~\delta +  \delta + \delta + \frac{\epsilon}{5} + \frac{\epsilon}{5} \\
			& \leq ~\frac{\epsilon}{5}  +\frac{\epsilon}{5} + \frac{\epsilon}{5}  + 2\frac{\epsilon}{5} 
			= \epsilon 
		\end{align*}
		Hence the Cauchy sequence $\{ x_n\}$ converges to $x\in X$, and this proves that $X$ is complete.
	\end{proof}
	%%%%%%%%%%%%%%%
	Cantor's intersection theorem ensures the completeness of a metric space.
	Our next theorem is the generalization of such theorem in a $\phi$-metric space.
	%%%%%%%%%%%%%%%%%%%%%
	\begin{thm}
		A necessary and sufficient condition that the $\phi$-metric space $ (X,d_\phi) $ be complete is that every nested sequence of non-empty closed subsets $ \{F_i\} $ with $ \delta(F_i)\rightarrow 0 $ as $i\rightarrow \infty $ be such that $F=\cap _{i=1} ^\infty F_i$ contains exactly one point. 
	\end{thm}
	\begin{proof}
		First suppose that $X$ is complete. Consider a sequence of closed subsets $\{F_i\}$ such that $ F_1\supset F_2\supset F_3\supset \dotsm  $ and $ \delta(F_i) \rightarrow 0 $ as $i\rightarrow\infty$. \\
		For all $ n\in \mathbb{N} $, choose $ a_n\in F_n$. Hence we generate a sequence $\{a_n\}$ in $X$. We verify that the sequence $ \{a_n\} $ is a Cauchy sequence.\\
		Now  for some $n\in \mathbb{N},~ a_n\in F_n$ implies $a_{n+p} \in F_{n+p}\subset F_n,~ \text{for all } ~ p=1,2, \cdots$. \\
		So, for all $ p=1,2, \cdots $,
		$ d_\phi(a_n,a_{n+p})\leq \delta(F_n), ~ \text{for all} ~ n\in \mathbb{N} ~ \text{which implies} ~ \underset{n\rightarrow\infty}{lim} d_\phi(a_n,a_{n+p})=0. $ 
		Hence,  $ \{a_n\} $ is a Cauchy sequence in $ X $.
		Since $X$ is complete, there exists $ a\in X$ such that \\
		$ a_n\rightarrow a $ as $n\rightarrow \infty$.\\
		For a fixed positive integer $k$,  consider the subset $F_k$. Then each $a_k,a_{k+1},a_{k+2},\dotsm  \in F_k$. Since $F_k$ is closed, so $a\in F_k$. $k$ being arbitrary positive integer, we can conclude $a\in \underset{i \in \mathbb{N}}{\cap}  F_i$.\\
		Finally we show that $a $ is unique. For, let there exists $ b(\neq a)\in \underset{i \in \mathbb{N}}{\cap}  F_i $. Then for each $k\in \mathbb{N}$, 
		$$ a,b\in F_k 
		~~ that ~ is ~~  d_\phi(a,b)\leq \delta(F_k) $$
		Therefore, $ d_\phi(a,b)=0,~~\text{since} ~\delta(F_k) \rightarrow 0~ \text{as} ~k\rightarrow \infty $, and hence $ a=b $. \\
		Conversely suppose that the condition of the theorem holds. 
		To show that $X$ is complete, take a Cauchy sequence $\{x_n\}$ in $ X$. \\
		Let $F_n=\{x_n,x_{n+1},x_{n+2}, \dotsm \} $, for all $ n\in \mathbb{N}$. If we choose any $\epsilon>0$, then there exists a positive integer $n_0$(say) such that 
		\begin{align*}
			&	d_\phi(x_n,x_m)<\epsilon, ~~~\text{for all} ~ n > m \geq n_0 \\
			or ~~ & \delta(F_n) \leq \epsilon, ~~\text{for all} ~ n \geq n_0  \\
			or ~~ & \delta(\overline {F_n})  \leq \epsilon, ~~\text{for all} ~ n \geq n_0  \\
			or ~~ & \delta(\overline {F_n}) \rightarrow 0 ~ ~\text{as}~ n\rightarrow \infty. 
		\end{align*}
		Clearly ${F_{n+1}}\subset {F_n} $, for each $n$, and thus  $ \overline {F_{n+1}} \subset \overline {F_n} $, for each $n$. So $\{ \overline {F_n} \}$ constitutes a closed, nested sequence of non-empty sets in $X$ whose diameter tends to zero. By hypothesis, there exists a unique point $x\in \underset{n \in \mathbb{N}}{\cap} \overline {F_n}$.\\
		Now for each $n=1,2,\dotsm, ~ x_n \in F_n \subseteq  \overline {F_n}$ implies
		$$	 d_\phi(x_n,x) \leq \delta(\overline{F_n}), ~~\text{for all} ~ n\in \mathbb{N}  $$
		Therefore, $ d_\phi(x_n,x)  \rightarrow 0 ~~ \text{as}~ n\rightarrow \infty $.  \\
		This shows that  the Cauchy sequence $ \{x_n\} $ converges to $ x \in X$, and hence $ X$ is complete.
	\end{proof} 
	%%
	%%%%%
	%%%%%%%
	%%%%%%%%%%   totally boundedness  **%%%%%%%%%%%%%
	Next, we define a property which is more potent than boundedness. 
	%%%%%%%%%%%%%%%%%%%%%%%%
	\begin{dfn}
		Let $(X,d_\phi)$ be a $\phi$-metric space and $B\subseteq X$. 	
		$B$ is said to be totally bounded if for every $\epsilon>0$, there exists a finite subset $A_\epsilon$ of $X$ such that $B=\underset{a\in A_\epsilon}{\cup} B_\phi(a,\epsilon)$.
	\end{dfn}
	%%%%%%%%%%%%%%%%%%%%%%%%%%%%%%%%%%
	\begin{thm}\label{tot bdd is bdd}
		Every  totally bounded subset of a $\phi$-metric space $(X,d_\phi)$ is bounded.  
		\begin{proof}
			Let $ \epsilon >0 $, and $B$ be a totally bounded subset of $(X,d_\phi)$.\\
			For any $\alpha, \beta \in B$, and for any $x,y\in X $, 
			\begin{align*}
				d_\phi (\alpha, \beta) & \leq d_\phi (\alpha, x) + d_\phi (x, \beta) + \phi (\alpha, \beta, x) \\
				&  \leq d_\phi (\alpha, x) + d_\phi (x, y) + d_\phi (y, \beta) + \phi( x, \beta, y) + \phi (\alpha, \beta, x) 
			\end{align*} 
			Choose	$ \epsilon = 1$.\\ 
			Then for $\epsilon =1 $, there exists $ \delta_1, \delta_2 >0$ such that for any $\alpha, \beta \in B$,  $ \phi( x, \beta, y) < 1 $ whenever \\
			$ d_\phi(y,\beta)< \delta_1,  ~x,y \in X $, and  $ \phi( \alpha, \beta, x) < 1$ whenever   $ d_\phi(x,\alpha)< \delta_2, ~ x\in X $.\\
			Let $\delta= \min \{ 1,  \delta_1, \delta_2\}$.\\
			Since  $B$ be a totally bounded, so for $\delta >0$, there exists a finite subset $ S=\{x_1,x_2,\cdots,x_n\} $ of $X$ such that $B=\cup _{x_i=1} ^n B_\phi(x_i,\delta)$.\\
			Choose any $\alpha, \beta \in B$.\\
			Then there exists $ x_i, x_j\in S $ such that $ \alpha \in B_\phi(x_i,\delta ) $, and $ \beta \in B_\phi(x_j,\delta) $.\\
			Hence we obtain, $  d_\phi (\alpha, x_i) < \delta,~  d_\phi (x_j, \beta) < \delta $, and $ \phi( x_i, \beta, x_j) < 1,~ \phi (\alpha, \beta, x_i) <1 $.\\
			Suppose $ K= \max \{ d_\phi(x_i,x_j): x_i, x_j \in S \} $. \\
			Therefore for any $ x_i, x_j \in S$ and $\alpha, \beta \in B$, 
			\begin{align*}
				d_\phi (\alpha, \beta) < & ~ \delta + d_\phi(x_i,x_j) + \delta + 1 +1 \\
				\leq & ~ 1 + \max \{ d_\phi(x_i,x_j): x_i, x_j \in S \}  + 1 + 2 
				=  ~ 4 + K
			\end{align*} 
			Since   $\alpha, \beta \in B$ is arbitrary, so   $  d_\phi (\alpha, \beta) <\leq K +4 $, for all $ \alpha, \beta \in B $. This proves that $B$ is bounded.
		\end{proof}	
	\end{thm}
	%%%%%%%%%%%%%%%%%%%%%
	\begin{thm}
		In a totally bounded $\phi$-metric space $(X,d_\phi)$,  every sequence has a Cauchy subsequence.
		\begin{proof}
			Let $\{ x_n \}$ be a sequence in $X$. Since $X$ is totally bounded so it can be covered by a finite number of open balls of any radius. Let us consider balls of radius $1$. Then at least one of these open balls, say $A_1$ contains infinitely many elements of the sequence. Choose $x_{k_1} \in A_1$  for some $k_1 \in \mathbb{N}$. Similarly, $A_1$ being totally bounded can be covered by a finite number of open balls each of  radius $ \frac{1}{2} $. Then at least one of these open balls, say $A_2$ contains infinitely many elements of the sequence. We choose $x_{k_2} \in A_2 $  for some $k_1<k_2 \in \mathbb{N}$. Continuing in this way we obtain a sequence $\{ A_n \}$ of open balls with radius $\frac{1}{n}$ such that $ A_1\supset A_2 \supset \cdots $ and $x_{k_n} \in A_{k_n} $ with  $k_1<k_2 \cdots $. Clearly $\{ x_{k_n} \}$ is a subsequence of $\{x_n\}$. Choose $\epsilon >0$. Then there exist $N \in \mathbb{N}$ such that $\frac{2}{N}< \epsilon$. Now for all $r,s\geq N, ~x_{k_r}, x_{k_s} \in A_N$, and hence $d_\phi (x_{k_r}, x_{k_s} )< \frac{2}{N}< \epsilon$ which implies $\{ x_{k_n} \}$ is a Cauchy sequence in $X$.
		\end{proof}	
	\end{thm}
	%%%%%%%%%%%%%%%%%%%%%%%%%%
	\begin{thm}\label{com is bdd}
		Every compact $\phi$-metric space $(X,d_\phi)$ is totally bounded.
		\begin{proof}
			From the compactness of $X$ it follows that, for every $\epsilon>0,~\beta = \{ B_\phi(a,\epsilon):a\in X \}$ is an open cover of $X$, and there exists a finite subset of $\beta$ which covers $X$. Therefore $X$ is totally bounded. 	
		\end{proof}
	\end{thm}
	%%%%%%%%%%%%%%%%%%%%%
	\begin{rem}
		The converse of the Theorem \ref{tot bdd is bdd} and Theorem \ref{com is bdd} does not hold in general. Since  metric spaces are also $\phi$-metric spaces, so this can be justified by the examples of a metric space.
	\end{rem}
	%%%%%%%%%%%%%%%%
	But totally boundedness and completeness together force the $\phi$-metric space to be compact. We prove this in our next theorem.
	%%%%%%%%%%%%%%%%%%
	\begin{thm}
		Let $(X,d_\phi)$ be a $\phi$-metric space. If $X$ is totally bounded and complete then $X$ is compact.
		\begin{proof}
			Let $\{ x_n \}$ be a Cauchy sequence in $X$. So totally boundedness of $X$ implies that $\{ x_n \}$ has a Cauchy subsequence, say $\{ x_{k_n} \}$. Since $X$ is complete, thus $\{ x_{k_n} \}$ converges in $X$. Therefore $X$ is sequentially compact, and hence compact.
		\end{proof}
	\end{thm}
	%
	%%%
	%%%%%%%%%%%%%%%%%%%%%%%%%%%%%%%%%
	%%%%%%%%%%%%%%%%%%%%%%%%%% Fixed point theories%%%%%%%%%%%%%%%%%
	%%%%%%%%%%%%%%%%%%%%%%%%%%%
	\section{Some fixed point theorems in $\phi$-metric spaces } 
	In this section, we establish the existence of a fixed point for the Banach type, and the Kannan type contraction principle, and also develop the Edelstein theorem in $\phi$-metric spaces.  \\
	Before going to the main results, we prove a useful lemma.  
	\begin{lma} \label{lmaforfxdpt}
		Let $(X,d_\phi)$ be a $\phi$-metric space. If $\{x_n\}_{n=0} ^\infty $ be a sequence in $X$ which satisfies 
		$$ d_\phi(x_n, x_{n+1}) \leq kd_\phi(x_{n-1},x_n), ~ n = 1, 2, \cdots $$
		where $0<k<1 $, then $\{x_n\}$ is a Cauchy sequence in $X$. 
		\begin{proof}
			Suppose $\{x_n\}_{n=0} ^\infty $ be a sequence in $X$ which satisfies the mentioned conditions. \\
			Let $ \epsilon >0$.\\
			Now for all $ n = 1, 2, \cdots $, 
			$$ d_\phi(x_n, x_{n+1}) \leq kd_\phi(x_{n-1},x_n)  $$ 
			implies 
			\begin{align*} 
				d_\phi(x_n, x_{n+1}) & \leq k^2 d_\phi(x_{n-2}, x_{n-1}) \\
				& \leq  k^3 d_\phi(x_{n-3}, x_{n-2}) \leq \cdots \leq  k^n d_\phi(x_0, x_1) 
			\end{align*} 
			Therefore, $ \underset{n\rightarrow\infty}{lim} d_\phi(x_n,x_{n+1})=0 $, since $ ~0<k<1  $.\\
			Again for $m>n$ we have, 
			\begin{align*}
				d_\phi(x_m,x_n)\leq & ~d_\phi(x_m,x_{n+1}) + d_\phi(x_{n+1},x_n) + \phi(x_m,x_n,x_{n+1})  \\
				\leq & ~\{d_\phi(x_m, x_{n+2}) + d_\phi( x_{n+2}, x_{n+1} ) + \phi (x_m, x_{n+1}, x_{n+2}) \}+ d_\phi( x_n, x_{n+1} ) + \\
				& \phi(x_m, x_n, x_{n+1})  \\
				\leq & ~ \{ d_\phi(x_m,x_{m-1})+...+ d_\phi( x_{n+1}, x_{n+2} )+ d_\phi( x_{n+1}, x_n) \} +\\
				& ~ \{ \phi(x_m, x_n, x_{n+1})+ \phi (x_m, x_{n+1}, x_{n+2}) +...+ \phi(x_m, x_{m-2}, x_{m-1}) \} \\
				\leq & ~ \{k^{m-1}+k^{m-2}+...+k^{n+1}+k^n\}d_\phi(x_1,x_0) + \sum_{i=n}^{{m-2}}\phi(x_m,x_i,x_{i+1}) \\
				\leq & ~ k^n\{1+ k+k^2+ \cdots \}d_\phi(x_1,x_0) + \sum_{i=n}^{{m-2}}\phi(x_m,x_i,x_{i+1}) 
			\end{align*} 
			which gives 
			\begin{equation} \label{contrctn seq}
				d_\phi(x_m,x_n) \leq \frac{k^n}{1-k} d_\phi( x_1, x_0 ) + \sum_{i=n}^{{m-2}} \phi(x_m,x_i,x_{i+1}), ~~  \text{for all} ~ m > n
			\end{equation}
			Now for $ \epsilon>0,$ there exists  $ \delta>0$ such that $ \phi(x_t, x_l, x_{l+1})<\epsilon$ whenever $d_\phi(x_l,x_{l+1})<\delta$. \\
			Let $ \beta = \min \{ \epsilon, \delta\} $. \\
			Then since $\underset{i\rightarrow\infty}{lim} d_\phi(x_i,x_{i+1})=0$, so for that $\beta>0 $, there exists a positive integer $ N $ such that $ d_\phi(x_i,x_{i+1})<\beta$, for all $ i \geq N $. \\ 
			Thus for all $ m > n \geq N, $ we have $ d_\phi(x_n,x_{n+1})<\beta $  and  $ \phi ( x_m, x_n, x_{n+1} ) < \epsilon $.  \\
			Since $ \epsilon>0$ arbitrarily chosen,  so $\underset{n\rightarrow\infty}{lim} d_\phi(x_n,x_{n+1})=0$ implies  $\underset{m,n\rightarrow\infty}{lim} \phi(x_m,x_n,x_{n+1})=0$. \\
			Hence the relation ( \ref{contrctn seq} ) gives $\underset{m, n\rightarrow\infty}{lim} d_\phi(x_n,x_m)=0 $
			which implies $\{x_n\}$ is a Cauchy sequence in $X$. 
		\end{proof}	
	\end{lma}
	%%%%%%%%%%%%%%%%%%%%%%%%%%%%%%%%
	\begin{thm} \label{banach}(Banach type contraction)
		Let $(X,d_\phi)$ be a complete $\phi$-metric space and $T$ be a self-mapping on $X$ satisfying 
		$$ d_\phi(Tx,Ty)\leq kd_\phi(x,y) $$ 
		for all  $x,y\in X $, where $k\in(0,1)$. Then $T$ has a unique fixed point in $X$.
		\begin{proof}
			If any fixed point of $T$ exists then uniqueness directly follows from the contraction condition. Here we only prove the existence of a fixed point. For,  
			consider an iterative sequence: \\  $x_0,~x_1=Tx_0,~x_2=Tx_1,...,x_{n+1}=Tx_n, \cdots $, for some fixed $x_0 \in X$.\\
			Then,
			$$ d_\phi(x_{n+1},x_n) = d_\phi(Tx_{n}, T {x_n-1 } ) \leq k	d_\phi(x_n, x_{n-1}), ~\text{for all} ~ n = 1, 2, \cdots . $$ 
			Hence by Lemma \ref{lmaforfxdpt}, we can conclude $\{x_n\}$ is a Cauchy sequence in $(X,d_\phi)$. Since  $X$ is complete, so $\{x_n\}$ converges to some $ \zeta \in X $. \\
			Lastly, we will prove that  $ \zeta $ is a fixed point for $T$. Now
			$$
			d_\phi(T\zeta,\zeta)= \underset{n\rightarrow\infty}{lim} d_\phi(T\zeta,x_n)\leq k \underset{n\rightarrow\infty}{lim} d_\phi(\zeta,x_{n-1})
			$$ 
			which implies $ d_\phi(T\zeta,\zeta) = 0 $ that is $ T\zeta=\zeta $.
		\end{proof}	
	\end{thm}
	%%%%
	%%
	%%
	\begin{eg}
		Consider the complete $\phi$-metric space $( X, d_\phi)$ where $X= \mathbb{R} $ and the $\phi$-metric defined by
		$ d_\phi(x,y)= (x-y)^2 $, for all $ x,y \in \mathbb{R} $.\\
		Let us define a self-mapping $T$ on $ \mathbb{R} $ defined by $ T x= \frac{\sin x}{2} $, for all $ x \in \mathbb{R} $.\\
		Then 
		$$ d_\phi(Tx, Ty)= \frac{1}{4} (\sin x-\sin y)^2 = [\cos(\frac{x+y}{2}) \cdot \sin (\frac{x-y}{2})]^2 \leq \frac{1}{4} ( x-y)^2 = \frac{1}{4} d_\phi(x,y) $$ 
		Thus $T$ satisfies the Banach type contraction for $ \frac{1}{4} \leq k <1 $, and $ x=0 $ is the unique fixed point for $ T $.
	\end{eg}
	%%
	%%%%
	%%%%
	%%%
	\begin{thm}\label{kannan}
		(Kannan type contraction)
		Let $(X,d_\phi)$ be a complete $\phi$-metric space and $T$ be a self-mapping on $X$ satisfying 
		$$ d_\phi(Tx,Ty)\leq k[d_\phi(x,Tx)+d_\phi(y,Ty)] $$ 
		for all  $x,y\in X $, where $k\in(0,\frac{1}{2})$.  Then $T$ has a unique fixed point in $X$.
	\end{thm}
	\begin{proof}
		For some fixed $x_0 \in X$,	consider the iterative sequence: \\  $x_0,~x_1=Tx_0,~x_2=Tx_1, \cdots ~,
		x_{n+1}=Tx_n, \cdots $. Then, 
		\begin{align*}
			& d_\phi(x_2,x_1)=d_\phi(Tx_1,Tx_0)\leq k[d_\phi(x_1,x_2)+d_\phi(x_0,x_1)]\\
			or ~~ & d_\phi(x_2,x_1)\leq \alpha d_\phi(x_1,x_0)~~\text{where}~~\alpha=\frac{k}{1-k}, ~\text{and}~0<\alpha<1
		\end{align*}
		Proceeding in this way, we can write  $d_\phi(x_{n+1},x_n)\leq \alpha d_\phi( x_n, x_{n-1} ) $, for all $ n \in \mathbb{N}$, where $ 0< \alpha<1 $. \\
		Applying Lemma  \ref{lmaforfxdpt}, we can conclude $\{x_n\}$ is a Cauchy sequence in $(X,d_\phi)$, and since  $X$ is complete, so $\{x_n\}$ converges to some $\zeta\in X$. Now, 
		\begin{align*}
			& d_\phi(T\zeta,\zeta)= \underset{n\rightarrow\infty}{lim} d_\phi(T\zeta,x_n)\leq k\underset{n\rightarrow\infty}{lim}[d_\phi(\zeta,T\zeta)+d_\phi(x_{n-1},x_n)]\\
			or ~~ & d_\phi(T\zeta,\zeta)\leq kd_\phi(T\zeta,\zeta)
		\end{align*}
		Therefore, $ d_\phi(T\zeta,\zeta)=0,~\text{as}~0<k<\frac{1}{2} $, which implies $ T\zeta=\zeta $. \\
		Hence $\zeta$ is a fixed point of $T$, and uniqueness easily follows from the contraction condition.  
	\end{proof}
	%%%%%%%%%%
	\begin{eg}
		Define a function $d_\phi$ on a set $ X= [-2,2] $ by $ d_\phi(x,y)= (x-y)^2, ~\forall x,y \in [-2,2] $. Then  $(X,d_\phi)$ is a complete  $\phi$-metric space.\\
		Define $ T: X \rightarrow X $ by 
		$$ T(x)=\begin{cases}
			\frac{x}{10} &\text{when}~ -2\leq x<1\\
			\frac{x}{5} &\text{when}~ ~~~ 1\leq x \leq 2
		\end{cases} $$
		For all $ x, y \in \mathbb{R}$ we know, 
		\begin{equation*}
			(|x|-|y|)^2 \geq 0  ~~ implies  ~~ |x|^2 + |y|^2  \geq  2|x| |y|
		\end{equation*}
		Hence we obtain, 
		\begin{align*}
			|x-y|^2 & \leq [|x| + |y|]^2 \\
			& = |x|^2 + |y|^2 + 2|x||y| \\
			& \leq (|x|^2 + |y|^2 ) + (|x|^2 + |y|^2),~\forall ~  x, y \in \mathbb{R}
		\end{align*}
		that is
		\begin{equation}\label{xxx}
			|x-y|^2  \leq 2(|x|^2 + |y|^2 ), ~ \text{for all} ~  x, y \in \mathbb{R}
		\end{equation}
		Case:I $~~$ Let $x, y \in [-2,1)$. Then 
		\begin{equation*}
			d_\phi(x,Tx)+d_\phi(y,Ty) = (x-\frac{x}{10})^2 + (y-\frac{y}{10})^2 = \frac{81}{100} [|x|^2 + |y|^2]
		\end{equation*}
		and 
		\begin{align*}
			d_\phi(Tx,Ty) & =\frac{1}{100} (x-y)^2 \\
			& \leq  \frac{2}{100} [|x|^2 + |y|^2 ], ~\text{using the relation }~(\ref{xxx})\\ 
			& = \frac{2}{81} [ d_\phi(x,Tx)+d_\phi(y,Ty) ]
		\end{align*}
		Case:II $~~$ Let $x, y \in [1,2]$. Then
		\begin{equation*}
			d_\phi(x,Tx)+d_\phi(y,Ty)   = (x-\frac{x}{5})^2 + (y-\frac{y}{5})^2  = \frac{16}{25} [|x|^2 + |y|^2]
		\end{equation*}
		and 
		\begin{align*}
			d_\phi(Tx,Ty) & =\frac{1}{25} (x-y)^2 \\
			& \leq  \frac{2}{25} [|x|^2 + |y|^2 ],  ~\text{using the relation }~(\ref{xxx})\\
			& = \frac{1}{8} [ d_\phi(x,Tx)+d_\phi(y,Ty) ]
		\end{align*}
		Case:III $~~$ Let $x\in [-2,1),~ y \in [1,2]$. Then
		\begin{equation*}
			d_\phi(x,Tx)+d_\phi(y,Ty)   = (x-\frac{x}{10})^2 + (y-\frac{y}{5})^2  = \frac{81}{100} |x|^2 + \frac{16}{25}  |y|^2
		\end{equation*}
		and 
		\begin{align*}
			d_\phi(Tx,Ty) & = (\frac{x}{10}-\frac{y}{5})^2 \\ 
			& \leq  2 [\frac{|x|^2}{100} + \frac{|y|^2}{25} ],  ~\text{using the relation }~(\ref{xxx}) \\
			& =  \frac{2}{81}\cdot \frac{81}{100} |x|^2 +\frac{2}{16} \cdot \frac{16}{25} |y|^2  \\
			& < \frac{1}{8} [ d_\phi(x,Tx)+d_\phi(y,Ty) ]
		\end{align*}
		Therefore for all $x, y \in [-2,2], ~ d_\phi(Tx,Ty)\leq k[d_\phi(x,Tx)+d_\phi(y,Ty)], $ where $ k= \frac{1}{8} $. 
	\end{eg}
	%%
	%%%%
	%%%%
	%%
	%%
	\begin{thm}\label{edelstien}
		Let $(X,d_\phi)$ be a $\phi$-metric space and $T$ be a self-mapping on $X$ satisfying
		$$ d_\phi(Tx,Ty)<d_\phi(x,y) $$
		for all  $x,y\in X $. If there exists  $ x \in X $ such that the sequence $ \{T^nx\} $ has a subsequence converging to $\zeta$ then $\zeta$ is the 
		unique fixed point of $T$.
		\begin{proof}
			Let $ \{T^{n_i}x\} $ be a subsequence of $ \{T^nx\} $ converging to $\zeta$. If $ T^kx = T^{k+1}x $ for some $k\in \mathbb{N}$,
			then $\zeta$ is a fixed point of $T$. \\
			If $ T^kx \neq T^{k+1}x $ for any $k\in \mathbb{N}$ and $T\zeta\neq \zeta$, then  
			\begin{equation}\label{eq9}
				d_\phi(T\zeta,T^2\zeta)<d_\phi(\zeta,T\zeta)
			\end{equation}
			For $x\in X$ and a fixed $ n_l$, for all $ n > n_l +1$ we have,
			\begin{equation*}
				d_\phi(T^nx,T^{n+1}x)<d_\phi(T^{n_l+1}x,T^{n_l+2}x)
			\end{equation*}
			Clearly $d_\phi(\zeta,T\zeta)$ is a limit point of the sequence $\{d_\phi(T^nx,T^{n+1}x)\}$ and so 
			\begin{equation}\label{eq10}
				d_\phi(\zeta,T\zeta)\leq d_\phi(T^{n_l+1}x,T^{n_l+2}x)
			\end{equation}
			In equation (\ref{eq10}) letting $l\rightarrow\infty$, we get $d_\phi(\zeta,T\zeta)\leq d_\phi(T\zeta,T^2\zeta)$, which contradicts equation
			(\ref{eq9}). So either $ T^kx = T^{k+1}x $ for some $k\in \mathbb{N}$ or $T\zeta= \zeta$ or both. Hence $\zeta$ is a fixed point of $T$, which is  unique also. 
		\end{proof}	
	\end{thm}
	%
	%%%
	%%%%%%
	%%%%
	\begin{eg}
		Consider the  $\phi$-metric space $( X, d_\phi)$ where $X= [-\frac{\pi}{2}, \frac{\pi}{2}] $ and the $\phi$-metric defined by
		$ d_\phi(x,y)= (x-y)^2$, for all $ x,y \in [-\frac{\pi}{2}, \frac{\pi}{2}] $.\\
		Let us define a self-mapping $T$ on $ [-\frac{\pi}{2}, \frac{\pi}{2}]$ defined by $ T x= \tan^{-1} x -x,$ for all $ x \in [-\frac{\pi}{2},
		\frac{\pi}{2}] $. Then,  
		\begin{align*}
			d_\phi(Tx, Ty)= & | \tan^{-1} x -x - \tan^{-1} y + y|^2 \\ 
			= & |(y-x)-(\tan^{-1}y-\tan^{-1}x) |^2 \\
			= & | (y-x)-(\frac{y-x}{1+ \eta^2}) |^2, ~\min \{x, y \} < \eta < \max \{x, y\} ~\text{ ( by Lagrange's mean value theorem ) }\\
			= & |y-x|^2\frac{\eta^2}{1+ \eta^2} \\
			< & ( x-y)^2 = d_\phi(x,y)
		\end{align*}
		Thus $T$ satisfies  contractive condition	 and for $ x=0 $, all subsequence of $ \{ T^n0\}=\{0\}$ converges to $0$. Hence from 
		theorem we conclude that $x=0$ is the unique fixed point for $T$.\\
		Moreover, from the defined self-mapping it is clear that $x=0$ is the unique fixed point for $ T $.	
	\end{eg}
	%%%%%%%%%%%%%%%%%%%%%%%%%
	\begin{rem}
		\begin{enumerate}[(a)]
			\item If we take a S-metric space $(X,S)$ then $d_\phi(x,y)=S(x,x,y)$, for all  $ x, y \in X $, defines a $\phi$-metric on $X$. Then
			Theorem \ref{banach} and Theorem \ref{edelstien} reduces to Theorem 3.1 and Theorem 3.3  of Section 3 of  Sedghi et al.\cite{[1]}. 
			\item If  we consider a b-metric space $(X,B)$ then $d_\phi(x,y)=B(x,y)$, for all $ x,y \in X $, defines a $\phi$-metric on $X$.
			Then Theorem \ref{banach} and Theorem \ref{kannan} reduce to Theorem 1 and Theorem 2  of the main results of   Mehmet Kir et al.\cite{[16]}. 
		\end{enumerate}
	\end{rem}
	~\\
	~\\
	%
	%
	%%%%%%%%%%%%
	\textbf{Conclusion:} This generalized metric, $ \phi$-metric is not only for seeking generalization. $\phi$–metric is helpful to study some
	existing generalized distance functions and may play the role of metrics in many scenarios. As we have mentioned earlier, it is very difficult to 
	find a metric for S-metric and b-metric spaces. The same fact is true for 
	$\phi$-metric spaces. However, one can easily  construct a $\phi$-metric from S-metric and b-metric which solves all the purposes of a metric
	for  topological spaces. So one does not need to find a metric at all. Another important fact is that, in the geometry of surfaces, straight lines
	do not satisfy the `Euclidean' rules in general. The triangle inequality is not necessarily true for non-planar surfaces and hence straight 
	line distance   is not a metric there. But it is a $ \phi$-metric. We have a plan to study non-Euclidean geometry using this $\phi $-metric 
	in our subsequent work. \\
	%%%%
	%%%
	~\\
	\textbf{Data Availability:} No data are used in this research. \\
	~\\
	\textbf{Conflicts of Interest:} The authors declare no conflicts of interest. \\
	~\\
	\textbf{Authors' Contribution:} All the authors contributed equally to this research work.\\

	~\\
	%%%%%%%%%%%%%%%%%%
	\textbf{Acknowledgment:} The author AD is thankful to University 
	Grant Commission (UGC), New Delhi, India for awarding her senior research fellowship [Grant No.1221/(CSIRNETJUNE2019)]. The authors are thankful to the Editor-in-Chief, Editors, and Reviewers of the journal (IJNAA) for their valuable comments which helped us to revise the manuscript in the present form. We are grateful to the Department of Mathematics, Siksha-Bhavana, Visva-Bharati. \\

\end{document}